\documentclass[fleqn,centertags]{amsart}

\newtheorem{theorem}{Theorem}%[section]

\newtheorem{conjecture}{Conjecture}

\newenvironment{thm}{\medskip \noindent\normalsize {\bf Theorem.}\it}{\medskip}
\newenvironment{prop}{\medskip \noindent\normalsize {\bf Proposition.}\it}{\medskip}
\newenvironment{lemma}{\medskip \noindent\normalsize {\bf Lemma.}\it}{\medskip}
\newenvironment{cor}{\medskip \noindent\normalsize {\bf Corollary.}\it}{\medskip}
\newenvironment{remark}{\noindent\normalsize {\bf Remark.}}{\medskip}

\def\SL{\operatorname{SL}}
\def\GL{\operatorname{GL}}
\def\SO{\operatorname{SO}}
\def\M{\operatorname{M}}

\def\Z{{\mathbb Z}}
\def\Q{{\mathbb Q}}
\def\R{{\mathbb R}}
\def\C{{\mathbb C}}
\def\A{{\mathbb A}}
\def\d{{\rm dim (\mathrm{G})/2}}

\def\cO{{\mathcal O}}
\def\cP{{\mathcal P}}
\def\I{{\mathcal I}}
\def\D{\mathcal{D}}

\def\G{{\mathrm G}}
\def\Gb{\overline{\mathrm G}}
\def\Gqs{{\mathcal G}}
\def\B{{\mathcal B}}
\def\Pb{\overline{P}}
\def\E{{\mathcal E}}
\def\F{{\mathcal F}}

\def\H{{\mathrm H}}

\def\An{\mathrm{A}}
\def\Bn{\mathrm{B}}
\def\Cn{\mathrm{C}}
\def\Dn{\mathrm{D}}
\def\En{\mathrm{E}}
\def\Fn{\mathrm{F_4}}
\def\Gn{\mathrm{G_2}}

\newcommand{\y}{\mathbf{y}}
\newcommand{\bfA}{\mathbf{A}}

\newcommand{\disc}{\mathcal{D}_K}
\newcommand{\Tr}{\mathrm{Tr}_{K/\Q}}

\newcommand{\order}{\mathcal{O}}
\DeclareFontFamily{OT1}{rsfs}{}
\DeclareFontShape{OT1}{rsfs}{n}{it}{<-> rsfs10}{}
\DeclareMathAlphabet{\mathscr}{OT1}{rsfs}{n}{it}
\def\beq{\begin{displaymath}}
\def\eeq{\end{displaymath}}

\begin{document}
\title{Counting maximal arithmetic subgroups}
\author[M. Belolipetsky]{Mikhail Belolipetsky\\
\\
with an appendix by Jordan Ellenberg and Akshay Venkatesh}
%\thanks{}
\address{
M.~Belolipetsky
\newline%\phantom{iiii}
Department of Mathematical Sciences,
Durham University,
%South Road,
Durham DH1 3LE,
U.K.
\newline%\phantom{iiii}
Institute of Mathematics,
Koptyuga 4,
630090 Novosibirsk, Russia
}
\email{mbel@math.nsc.ru}

\address{
J.~Ellenberg
\newline%\phantom{iiii}
Department of Mathematics,
Princeton University,
% 808 Fine Hall,
% Washington Road,
Princeton, NJ 08544,
U.S.A.
}
\email{ellenber@math.princeton.edu}

\address{
A.~Venkatesh
\newline%\phantom{iiii}
Department of Mathematics,
Massachusetts Institute of Technology,
%77 Mass. Ave,
\newline%\phantom{iiii}
Cambridge, MA 02139-4307, U.S.A.
}
\email{akshayv@math.mit.edu}

%\date{\today}

\begin{abstract}
We study the growth rate of the number of maximal arithmetic subgroups of bounded covolumes
in a semi-simple Lie group using an extension of the method due to Borel and Prasad.
\end{abstract}

\maketitle

\section{Introduction}
A classical theorem of Wang~\cite{W} states that a simple Lie group not locally isomorphic to
$\SL_2(\R)$ or $\SL_2(\C)$ contains only finitely many conjugacy classes of discrete subgroups of
bounded covolumes. This theorem describes the distribution of lattices in the higher rank Lie groups,
it also brings an attention to the quantitative side of the distribution picture. By now several
attempts have been made towards a quantitative analogue of Wang's Theorem but still very little
is known.

The problem of determining the number of discrete subgroups of bounded covolumes naturally
splits into two parts: the first part is to count maximal lattices and the second,
to count subgroups of bounded index in a given lattice. A recent project initiated by
A.~Lubotzky in~\cite{L} resulted in a significant progress towards understanding the subgroup
growth of lattices, which also allowed to formulate general conjectures on the asymptotic of
the number of lattices of bounded covolumes in semi-simple Lie groups (see~\cite{BGLM}, \cite{GLP},
\cite{LN} and \cite{LS}). In this paper we consider
another part of the problem: we count maximal lattices in a given semi-simple Lie group.
\medskip

Let $H$ be a product of groups $H_s(k_s)$, $s\in S$, where $S$ is a finite set, each $k_s$ is an
archimedean local field (i.e. $k_s = \R$ or $\C$) and $H_s(k_s)$ is an absolutely almost simple
$k_s$-group. Then $H$ is a semi-simple Lie group. Throughout this paper we shall consider only
semi-simple Lie groups of this form. We shall assume, moreover, that $H$ is connected and none
of the factors $H_s(k_s)$ is compact or has type $\An_1$. In particular, $H$ can be a non-compact
simple Lie group (real or complex) not locally isomorphic to $\SL_2(\R)$ or $\SL_2(\C)$.

Let $m^u_H(x)$ and $m^{nu}_H(x)$ denote the number of conjugacy classes of maximal cocompact irreducible
arithmetic subgroups and maximal non-cocompact irreducible arithmetic subgroups in $H$ of covolume less
than $x$, respectively. If the real rank of $H$ is greater than $1$, then by Margulis' theorem
\cite[Theorem~1, p.~2]{M} these numbers are equal to the numbers of the conjugacy classes of maximal
uniform and non-uniform  irreducible lattices in $H$ of covolume less than $x$. In the real rank $1$
case there may also exist non-arithmetic lattices which we do not consider here.

\begin{theorem} {\bf A.} If $H$ contains an irreducible cocompact arithmetic subgroup (or, equivalently,
$H$ is isotypic), then there exist effectively computable positive constants $A$ and $B$ which
depend only on the type of almost simple factors of $H$, such that for sufficiently large $x$
$$
x^A \le m^u_H(x) \le x^{B\beta(x)},
$$
where $\beta(x)$ is a function which we define for an arbitrary $\epsilon > 0$ as
$\beta(x) = C(\log x)^\epsilon$, $C = C(\epsilon)$ is a constant which depends only on $\epsilon$.
\medskip

\noindent{\bf B.}
If $H$ contains a non-cocompact irreducible arithmetic subgroup then there exist effectively computable
positive constants $A'$, which depends only on the type of almost simple factors of $H$, and $B'$ depending
on $H$, such that for sufficiently large~$x$
$$
x^{A'} \le m^{nu}_H(x) \le x^{B'}.
$$
\end{theorem}
%\beta(x) = (\log x)^{0.51};\\

Conjecturally, the function $\beta(x)$ in part~A can be also replaced by a constant and the constant $B'$
in part~B depends only on the type of almost simple factors of $H$. This would require,
in particular, a polynomial bound on the number of fields with a bounded
discriminant. The existence of such a bound is an old conjecture in number theory which is
possibly due to Linnik, it appears in a stronger form as Conjecture~9.3.5 in
H.~Cohen's book~\cite{Coh}. In fact, as a corollary from the proof of Theorem~1
we can show an equivalence of the conjecture ``$\beta(x)= {\rm const}$'' to
Linnik's conjecture (see Sect.~\ref{s65}).
\medskip

Theorem~1 is motivated by the problem of distribution of lattices in semi-simple Lie groups.
An application of this results and its corollaries in Section~6 to the problem is a part
of a joint work with A.~Lubotzky~\cite{BL}.
\medskip

Concluding the introduction let us briefly outline the proof of the theorem.
Our method is based on the work of Borel and
Prasad~\cite{BP}. What makes our task different from the results in~\cite{BP}
is that besides proving the finiteness of the number of arithmetic subgroups of bounded covolumes
we want to give bounds or, at least, {\nobreak asym}ptotic bounds for the number. This requires
certain modifications to the method on one side and some special number theoretic results on
the other. In Proposition~\ref{prop:BP} we improve a number theoretic result from
\cite[Sect.~6]{BP} which enables us to effectively count the possible fields of definition
of the maximal arithmetic subgroups. The proof of this proposition is technical and can be
safely skipped in the first reading. The key ingredient for a good upper bound for the number
of fields is an extension of a recent work of J.~Ellenberg and A.~Venkatesh~\cite{EV}, which
we formulate in Proposition~\ref{prop:ND} and for which the proof is given in Appendix.
After bounding the number of possible fields of definition $k$, we count admissible
$k$-forms, corresponding collections of local factors and conjugacy classes of arithmetic
subgroups. Here we use some Galois cohomology techniques, the Hasse principle
and basic number theoretic Proposition~\ref{prop:NP}. Altogether these lead to the proof
of the upper bounds in Theorem~1 given in Section~4. The lower bounds are easier, they are
established in Section~5.
\medskip

\noindent
{\bf Acknowledgement.} I am grateful to Alex Lubotzky for his lecture, which inspired this
project, and for many helpful discussions. I wish to thank Tsachik Gelander, Alireza
Salehi Golsefidy,  Gopal Prasad and the anonymous referees for the critical remarks and suggestions.
Jordan Ellenberg and Akshay Venkatesh contributed to this project valuable number theoretic
support. This research was carried out while I held a Golda--Meir Postdoctoral Fellowship
at the Hebrew University of Jerusalem.
I also would like to thank for hospitality the MPIM in Bonn where a part of this work was
completed.

\section{Preliminaries on arithmetic subgroups}\setcounter{equation}{0}

The aim of this section is to give a short account of the fundamental results of
Borel and Prasad (\cite{BP}, \cite{P}) which will be used in this paper. We would like to
encourage the reader to look into the original articles cited above for a better understanding
of the subject. Our modest purpose here is to settle down the notation and recall some formulas
for the future reference.

\subsection{}
Throughout this paper $k$ is a number field, $\cO_k$ its ring of integers, $V = V(k)$ the set
of places (valuations) of $k$ which is the union of $V_\infty(k)$ archimedean and $V_f(k)$ nonarchimedean
places, and $\A = \A_k$ is the ring of ad\`eles of $k$. The number of archimedean places of $k$
is denoted by $a(k) = \#V_\infty(k)$. Let $r_1(k)$, $r_2(k)$ denote the number of real and
complex places of $k$, respectively, so $a(k) = r_1(k) + r_2(k)$.
As usual, $\D_k$ and $h_k$ stay for the absolute value of the discriminant of
$k /\Q$ and the class number of $k$. For a finite extension $l$ of $k$, $\D_{l/k}$ denotes the $\Q$-norm
of the relative discriminant of $l$ over $k$.
\medskip

All logarithms in the paper are in base $2$ unless stated otherwise.
%\medskip

% \subsection{}
% If a connected semi-simple Lie group have a finite center (which is always the case for the
% linear groups), then there is a correspondence between arithmetic subgroups of $H$ and of the
% adjoint group ${\overline H} = {\rm Ad\:} H$ (see~\cite[Chapter~9, Def.~6.3]{M}). This allows us
% to reduce the problem of counting arithmetic subgroups of $H$ to ${\overline H}$: although some subgroups
% which are not conjugate in $H$ may correspond to conjugate subgroups of ${\overline H}$ this will not affect
% our considerations. So from now on we shall assume that $H$ is a connected semi-simple Lie group with
% a trivial center. This implies in particular that $H$ can be identified with the identity component
% of an $\R$-group $\G$.

\subsection{}\label{s21} Let $\G/k$ be an algebraic group defined over a number field $k$
such that there exists a continuous surjective homomorphism $\phi: \G(k\otimes_\Q\R)^o\to H$
with a compact kernel. We shall call such fields $k$ and $k$-groups $\G$ {\it admissible}.
If $S\subset V_\infty(k)$ is the set of archimedean places of $k$ over
which $\G$ is isotropic (i.e., non-compact), then $\phi$ induces an epimorphism
$\G_S = \prod_{v\in S} \G(k_v)^o\to H$ whose kernel is a finite subgroup of $\G_S$.

We consider $\G$ as a $k$-subgroup of $\GL(n)$ for large enough $n$ and define a subgroup
$\Gamma$ of $\G(k)$ to be {\it arithmetic} if it is commensurable with the subgroup of $k$-integral
points $\G(k)\cap \GL(n,\cO_k)$, that is, the intersection $\Gamma\cap\GL(n,\cO_k)$ is of finite
index in both $\Gamma$ and $\G(k)\cap\GL(n,\cO_k)$. The subgroups of $H$ which are commensurable
with $\phi(\Gamma)$ for some admissible $\G/k$ are called {\it arithmetic subgroups} of $H$ defined
over the field $k$.

We restrict ourselves to the irreducible lattices which implies that in the definition of the arithmetic
subgroups it is enough to consider only simply connected absolutely almost simple algebraic groups
$\G$ (see~\cite[Chapter~9.1]{M}).

\subsection{}\label{s22} A semi-simple Lie group
contains irreducible lattices if and only if all its almost simple factors have the same type ($H$
is {\it isotypic}). For example, we can take $H = \SL(2,\R)^a\times\SL(2,\C)^b$ or
$H = \SO(p_1,q_1)\times\SO(p_2,q_2)$ ($p_1+q_1 = p_2+q_2$), but in $\SL(2,\R)\times\SL(3,\R)$ all
lattices are reducible. Sufficiency of this condition is provided by the Borel--Harder
theorem~\cite{BH} and its necessity is discussed e.g. in \cite[Chapter~9.4]{M}.
Note that in general the assumption that $H$ is isotypic does not imply that $H$ contains non-uniform
irreducible lattices as is shown by an example due to G.~Prasad (see Prop.~12.31 in D.~Witte Morris,
{\em Introduction to Arithmetic Groups}, preliminary version of the book available from the author's
website at http://people.uleth.ca/$\tilde{\ }$dave.morris). This is the
reason why we have to impose an additional assumption concerning the existence of non-uniform irreducible
lattices in part~B of the theorem.

\subsection{}\label{s24} The methods of Borel--Prasad require a considerable amount of the
Bruhat--Tits theory of reductive groups over local fields. We assume familiarity with the
theory and recall only some basic definitions. An extensive account of what we need can be
found in Tits' survey article~\cite{T}.

Let $K$ be a nonarchimedean local field of characteristic zero (a finite
extension of the $p$-adic field $\Q_p$) and let $\G$ be an absolutely almost
simple simply connected $K$-group. The Bruhat--Tits theory associates to $\G/K$
a simplicial complex $\B = \B(\G/K)$ on which $\G(K)$ acts by simplicial
automorphisms which are special (which implies, in particular, that if an
element of $\G(K)$ leaves a simplex of $\B$ stable, then it fixes the simplex
pointwise). The complex $\B$ is called the {\it affine building} of $\G/K$. A
{\it parahoric subgroup} $P$ of $\G(K)$ is defined as a stabilizer of a simplex
of $\B$. Every parahoric subgroup is compact and open in $\G(K)$ in the
$p$-adic topology. Minimal parahoric subgroups are called {\it Iwahori}, they
are defined as subgroups of $\G(K)$ fixing chambers (i.e., maximal simplexes)
in $\B$. All Iwahori subgroups are conjugate in $\G(K)$. Maximal parahoric
subgroups are the maximal compact subgroups of $\G(K)$, they are characterised
by the property of being stabilizers of the vertices of $\B$. A maximal
parahoric subgroup is called {\it special} if it fixes a {\it special vertex}
of $\B$. A vertex $x\in\B$ is special if the affine Weyl group $W$ of $\G(K)$
is a semidirect product of the translation subgroup by the isotropy group $W_x$
of $x$ in $W$. In this case $W_x$ is canonically isomorphic to the (finite)
Weyl group of the $K$-root system of $\G$. If $\G$ is quasi-split over $K$ and
splits over an unramified extension of $K$, then $G(K)$ contains {\it
hyperspecial} parahoric subgroups. We refer to \cite[1.10]{T} for the
definition of hyperspecial parahorics, an important property of these subgroups
is that they have maximal volumes among all parahoric subgroups
\cite[3.8.2]{T}.

\subsection{}\label{s25}
We shall now define a Haar measure $\mu$ on $H$ with respect to which the volumes of arithmetic quotients
will be computed. Of course, the final result then will hold for any other normalization of the Haar
measure on $H$. The definition and most of the subsequent facts come from \cite{P}, \cite{BP}.

Let $\G$ be a semi-simple algebraic group defined over a number field $k$. Let $v\in V_\infty(k)$
be an archimedean place of $k$. First consider the case when $k_v = \R$. There exists a unique
anisotropic $\R$-form $\G_{cpt}$ of $\G$ which has a natural Haar measure giving the group
volume $1$. This measure can be transferred to $\G(k_v)$ in a standard way and we define
$\mu_v$ as its image. It is a canonical Haar measure on $\G(\R)$. In case $k_v = \C$ we have
$\G(k_v) = \G_1(\R)$ with $\G_1 = {\rm Res}_{\C/\R}\G$, and we define $\mu_v$ to be equal to
the canonical measure on $\G_1(\R)$. For a subset $S$ of $V_\infty(k)$ the Haar measure $\mu_S$
on $\G_S$ is defined as a product of $\mu_v$, $v\in S$. This also gives us a measure $\mu$ on
a semi-simple Lie group $H$.

\subsection{}\label{s26}
A collection $P = (P_v)_{v\in V_f}$ of parahoric subgroups $P_v$ of a simply connected $k$-group
$\G$ is called {\it coherent} if
$\prod_{v\in V_\infty}\G(k_v)\cdot\prod_{v\in V_f} P_v$ is an open subgroup of the ad\`ele group
$\G(\A_k)$. A coherent collection of parahoric subgroups $P = (P_v)_{v\in V_f}$ defines an arithmetic
subgroup $\Lambda = \G(k)\cap\prod_{v\in V_f} P_v$ of $\G(k)$ which will be called {\it
principal arithmetic subgroup} associated to $P$. The corresponding arithmetic subgroup
$\Lambda' = \phi(\Lambda)\subset H$ shall be also called {\it principal}.

The covolume of a principal arithmetic subgroup $\Lambda$ of $\G_S$ (in our case $S\subset V_\infty$)
with respect to the measure $\mu$ defined as above is given by Prasad's formula~\cite[Theorem~3.7]{P}:
$$ % \mu(H/\Lambda')  =
\mu_S(\G_S/\Lambda) = \D_k^{\d}(\D_l/\D_k^{[l:k]})^{\frac12s}
\left(\prod_{i=1}^{r}\frac{m_i!}{(2\pi)^{m_i+1}}\right)^{[k:\Q]} \tau_k(\G)\:\E(P),
$$
where
\begin{itemize}
\item[] ${\rm dim (\G)}$ and $m_i$ denote the dimension and Lie exponents of $\G$;
\smallskip\item[] $l$ is a Galois extension of k defined as in \cite[0.2]{P}
(if $\G$ is not a $k$-form of type $^6\Dn_4$, then $l$ is the split field of the quasi-split
inner $k$-form of $\G$, and if $\G$ is of type $^6\Dn_4$, then $l$ is a fixed cubic
extension of $k$ contained in the corresponding split field; in all the cases $[l:k]\le 3$);
\smallskip\item[] $s = s(\G)$  is an integer defined in  \cite[0.4]{P},
in particular, $s=0$  if  $\G$  is an inner form of a split group and $s\ge 5$ if $\G$ is an outer form;
\smallskip\item[] $\tau_k(\G)$ is the Tamagawa number of $\G$ over $k$, since $\G$ is simply connected
and $k$ is a number field $\tau_k(\G) = 1$;
\smallskip\item[]
$\E(P) = \prod_{v\in V_f} e_v$ is an Euler product of the local factors $e_v = e(P_v)$;
for $v\in V_f$, $e_v$ is the inverse of the volume of $P_v$ with respect to the Haar measure
$\gamma_v\omega_v^*$ defined in \cite[1.3, 2.1]{P}. These factors can be computed using the
Bruhat--Tits theory, in particular, $e_v > 1$ for every $v\in V_f$ \cite[Prop.~2.10(iv)]{P}.
\end{itemize}

\subsection{}\label{1:BP}
Any maximal arithmetic subgroup $\Gamma$ of $H$ can be obtained as a normalizer in $H$ of the
image $\Lambda'$ of some principal arithmetic subgroup of $\G(k)$ (see \cite[Prop.~1.4(iv)]{BP}).
Moreover, the collections of parahoric subgroups which are associated to the maximal arithmetic subgroups
have maximal types in a sense of Rohlfs~\cite{R} (see also~\cite{CR}). So in order to compute the
covolume of a maximal arithmetic subgroup we need to be able to compute the index
of a principal arithmetic subgroup in its normalizer. In a general setting the upper bound for the
index was obtained in~\cite[Sect.~2]{BP}:
\begin{equation*}
[\Gamma : \Lambda'] \le n^{\epsilon\# S}\cdot \#\H^1(k,{\rm C})_\xi \cdot \prod_{v\in V_f}\#\Xi_{\Theta_v}.
\end{equation*}
% The set $\S$ consists of all the archimedean places of $k$ over which $\G$ is not compact,
Here $n$ and $\epsilon$ are constants to be defined below, so $n^{\epsilon\# S}$ depends only
on $H$ and does not depend on the choice of $\G(k)$ and $\Lambda$. $\H^1(k,{\rm C})_\xi$ is a finite
subgroup of the first Galois cohomology group of $k$ with coefficients in the center of $\G$
defined in \cite[2.10]{BP}. The order of $\H^1(k,{\rm C})_\xi$ can be further estimated (see
\cite[Sect. 5]{BP}), the following bound is a combination of Propositions 5.1 and 5.6 [loc.~cit.]:
\begin{equation*}
\# \H^1(k,{\rm C})_\xi \le 2h_l^{\epsilon'} n^{\epsilon a(k) + \epsilon' a(l) +\epsilon \# T} (\D_l/\D_k^{[l:k]})^{\epsilon''},
\end{equation*}
where
\begin{itemize}
\item[]
$n = r+1$ if $\G$ is of type $\An_r$; $n = 2$ if $\G$ is of type $\Bn_r$, $\Cn_r$ ($r$ arbitrary),
$\Dn_r$ (with $r$ even), or $\En_7$; $n = 3$ if $\G$ is of type $\En_6$; $n = 4$ if $\G$ is
of type $\Dn_r$ with $r$ odd; $n = 1$ if $\G$ is of type $\En_8$, $\Fn$ or $\Gn$;
% BP:2.6\smallskip
\item[]
$\epsilon = 2$ if $\G$ is of type $\Dn_r$ with $r$ even, $\epsilon = 1$ otherwise
% BP:2.6
(so the center ${\rm C} = {\rm C}(\G)$ is isomorphic to $(\Z/n\Z)^\epsilon$ and $\#{\rm C} = n^{\epsilon}$);
\smallskip\item[]
$\epsilon' = \epsilon$ if $\G$ is an inner form of a $k$-split group, $\epsilon' = 1$ otherwise;
\smallskip\item[]
$\epsilon'' = 1$ if $\G/k$ is an outer form of type $\Dn_r$ ($r$ even) and $=0$ otherwise;
\smallskip\item[]
$T$ is the set of places $v\in V_f$ for which $\G$ splits over an unramified extension of $k_v$
but is not quasi-split over $k_v$.
\end{itemize}
Finally, $\Xi_{\Theta_v}$ is a subgroup of the automorphism group of the affine Dynkin diagram
which is coming from the adjoint group and preserving the type $\Theta_v$ of $P_v$. In particular,
$\#\Xi_{\Theta_v}\le r+1$ and $\#\Xi_{\Theta_v} = 1$ if $P_v$ is special.

\subsection{}\label{2:BP} As a result we have the following lower bound for the covolume of $\Gamma$:
$$
\mu(H/\Gamma) \ge (2n^{\epsilon a(k)+\epsilon' a(l)}h_l^{\epsilon'})^{-1}
\D_k^{\d}\left(\frac{\D_l}{\D_k^{[l:k]}}\right)^{\hspace*{-.4em}s'}\hspace*{-.6em}
\left(\prod_{i=1}^{r}\frac{m_i!}{(2\pi)^{m_i+1}}\right)^{[k:\Q]}\hspace*{-1em}\tau_k(\G)\F,
$$
where
\begin{itemize}
%BP:7.4
\item[]
$s' = s/2 - 1$ if $\G/k$ is an outer form of type $\Dn_r$, $r$ even, and $s' = s/2$ otherwise;
%BP:7.4
\smallskip\item[]
$\F = \prod_{v\in V_f} f_v$, with $f_v = e_v(\#\Xi_{\Theta_v})^{-1} = e_v$ if $\G$
is quasi-split over $k_v$ and $P_v$ is hyperspecial (which is true for almost all $v$),
$f_v = e_vn^{-\epsilon}(\#\Xi_{\Theta_v})^{-1}$ if $\G$ splits over an unramified extension of $k_v$
but is not quasi-split over $k_v$,
% $\G/k$ is an inner form of a split group and $\G$ does not split over $k_v$ or $\G/k$ is an outer
% form and $\G$ splits over an unramified extension of $k_v$ but is not quasi-split over $k_v$,
and $f_v = e_v(\#\Xi_{\Theta_v})^{-1}$ in the rest of the cases.
Using the computations in~\cite[Appendixes A,C]{BP} it is not hard to check that $f_v > 1$ for every
$v\in V_f$.
\end{itemize}
\medskip

\noindent
We refer to \cite[Sections 3,5]{BP} for more details about this formula.

% Note that our notation $f_v$ differ from the one used in the proof in \cite[7.4]{BP} by factor
% $(\#\Xi_{\Theta_v})^{-1}$ but still the argument in~\cite[Appendixes A,C]{BP} implies $f_v > 1$.

% Here we take a more precise and slightly more complicated estimate that the one which was mainly
% used in the proof of the finiteness theorems in~\cite{BP}. The reason for this is that we shall
% need some additional control over the volumes.

\section{Number theoretic results}\setcounter{equation}{0}

\subsection{}\label{prop:ND}
Let $N_{k,d}(x)$ be the number of $k$-isomorphism classes of extensions $l$ of
$k$ such that $[l:k] = d$, $\D_{l/k} < x$, and let $N(x)$ be the number of
isomorphism classes of number fields with discriminant less than $x$.

\begin{prop} For large enough positive $x$ we have
\begin{itemize}
\item[$(i)$] given a number field $k$ and a fixed degree $d$, there exist absolute constants
$c$, $b_1$, $b_2 > 0$ such that $N_{k,d}(x)\le c\D_k^{b_1}x^{b_2}$;
\item[$(ii)$] for every $\epsilon > 0$ there exists a constant $C = C(\epsilon) > 0$ such that \\
$N(x) \le x^{\beta(x)}$, $\beta(x) = C(\log x)^{\epsilon}$.
\end{itemize}
\end{prop}
\begin{proof} $(i)$ follows e.g. from \cite[Theorem~1.1]{EV}.
\medskip

\noindent
$(ii)$ follows from the method of~\cite{EV} but requires some extra work, namely, we have to
know how the implicit constants in \cite[Theorem~1.1]{EV} depend on the degree of the extensions
in order to be sure that this does not change the expected upper bound. This will be carried out
in detail in Appendix provided by J.~Ellenberg and A.~Venkatesh.
\end{proof}

\subsection{}\label{prop:NP}
Let $Q_k(x)$ be the number of  squarefree ideals of $k$ of norm $\le x$.
\medskip

\begin{prop}
\begin{itemize}
\item[$(i)$] $\displaystyle Q_k(x) = \frac{{\rm Res}_{s=1}(\zeta_k)}{\zeta_k(2)}\ x + o(x)$ for $x\to\infty$;
\item[$(ii)$] There exist absolute constants $b_3$, $b_4$ (not depending on $k$) such that\\
$Q_k(x) \le \D_k^{b_3}x^{b_4}$.
\end{itemize}
\end{prop}
\begin{proof}$(i)$ is a known fact from the analytic number theory. For a short and conceptual
proof we refer to \cite[Theorem 14]{S}.
\medskip

\noindent
$(ii)$ As far as we do not claim that $b_4 = 1$ the proof is easy. Consider the
Dedekind zeta function of $k$:
\begin{equation*} \zeta_k(s) = \sum_{n=1}^{\infty}\frac{a_n}{n^s},\end{equation*}
$a_n$ is the number of ideals of $k$ of norm $n$, $s>1$.

Let $I_k(x)$ denote the number of ideals of $k$ of norm less than $x$. We have
\begin{gather*}
I_k(x) = a_1 + a_2 + \ldots + a_{[x]};\\
\zeta_k(s)\cdot x^s \ge I_k(x).
\end{gather*}
Taking $s = 2$ we obtain
\begin{equation*}
Q_k(x) \le I_k(x) \le \zeta_k(2)\cdot x^2 \le \left(\frac{\pi^2}{6}\right)^{[k:\Q]} x^2 \le c_1^{\log \D_k} x^2
= \D_k^{b_3}x^{b_4}.
\end{equation*}
Here we used inequalities $\zeta_k(2)\le\zeta(2)^{[k:\Q]}$ and $[k:\Q] \le c\log \D_k$. The first
of them follows from the definition of the functions $\zeta$ and $\zeta_k$, and the second is a well
known corollary of Minkowski's discriminant bound.
\end{proof}

\subsection{} \label{prop:BP}
Finally we shall need an improved version of a number theoretic result from~\cite[Sect.~6]{BP}.
The main idea is that instead of using $\D_k^\d$ to absorb the small factors in the volume
formula we shall use only part of it saving the rest for a later occasion. This is easy to achieve
for the groups of a large enough absolute rank while when the rank becomes small the estimates become
much more delicate.

Let $\G/k$ be an absolutely almost simple, simply connected algebraic group of absolute rank $r\ge2$, so the
numbers $n$, $\epsilon$, $\epsilon'$, $s'$ and $m_1\le\ldots\le m_r$ are fixed and defined as in
Sect.~2. Let
\begin{align*}
B(\G/k) = \D_k^{\d} n^{-\epsilon a(k)-\epsilon' a(l)}h_l^{-\epsilon'}
(\D_l/\D_k^{[l:k]})^{s'}\left(\prod_{i=1}^{r}\frac{m_i!}{(2\pi)^{m_i+1}}\right)^{[k:\Q]}\!\!\!\!\!\!\!\!\!.
\end{align*}
Then by \ref{2:BP} we have $\mu(H/\Gamma) \ge \frac12 B(\G/k)\tau_k(\G)\F \ge \frac12 B(\G/k)$ for every
arithmetic subgroup $\Gamma$ of $H$ which is associated to $\G/k$.
\medskip

\begin{prop}
There exist positive constants $\delta_1$, $\delta_2$ depending only on the absolute type of $\G$,
such that $B(\G/k) \ge \D_k^{\delta_1} \D_{l/k}^{\delta_2}$ for almost all number fields $k$.
\end{prop}

\begin{proof}
Given an absolutely almost simple, simply connected algebraic group $\G$ of an absolute type $T$
and rank $r$, we shall show that for almost all $k$:
\begin{itemize}
\item[$(i)$] $B(\G/k) \ge \D_k^{\d-2} \D_{l/k}$ \quad if $r\ge 30$;
\item[$(ii)$] $B(\G/k) \ge \D_k \D_{l/k}$       \quad if $r<30$ and $T$ is not $\An_2,\ \An_3,\ \Bn_2$;
\item[$(iii)$] $B(\G/k) \ge \D_k^{0.1} \D_{l/k}$ \ if $T$ is $\An_3$ or $\Bn_2$;
\item[$(iv)$] $B(\G/k) \ge \D_k^{0.01} \D_{l/k}^{0.5}$  if $T$ is $\An_2$.
\end{itemize}
Clearly, altogether these four inequalities imply the proposition.
\medskip

Assume first that $\G$ is not a $k$-form of type $^6\Dn_4$. We have
\begin{gather*}
[l:k]\le 2;\\
n^{-\epsilon a(k)-\epsilon' a(l)} \ge n^{-\epsilon( a(k)+ a(l))} \ge (r+1)^{-3[k:\Q]}.
\end{gather*}
It is known that
\begin{equation}\label{31} h_l \le 10^2\left(\frac{\pi}{12}\right)^{[l:\Q]}\D_l \end{equation}
(cf.~\cite[proof of 6.1]{BP}, let us point out that this bound holds without
any assumption on the degree of the field $l$);
\begin{equation}\label{32} \D_l/\D_k^{[l:k]} = \D_{l/k} \ge 1.\end{equation}
%(follows e.g. from Theorem~A in the Appendix of~\cite{P}).

Combining the above inequalities we obtain
\begin{align*}
\lefteqn{
B(\G/k) = \D_k^{\d}n^{-\epsilon a(k)-\epsilon' a(l)}h_l^{-\epsilon'}
(\D_l/\D_k^{[l:k]})^{s'}\left(\prod_{i=1}^{r}\frac{m_i!}{(2\pi)^{m_i+1}}\right)^{[k:\Q]}} \\
 &\quad \ge 10^{-2\epsilon'}\D_k^{\d - 2}\left(\frac{\pi}{12}\right)^{-\epsilon'[l:\Q]}
   \D_{l/k}\left(\frac{1}{(r+1)^3}\prod_{i=1}^{r}\frac{m_i!}{(2\pi)^{m_i+1}}\right)^{[k:\Q]} \\
 &\quad \ge 10^{-2\epsilon'}\D_k^{\d - 2}\D_{l/k}
   \left(\frac{1}{(\pi/12)(r+1)^3}\prod_{i=1}^{r}\frac{m_i!}{(2\pi)^{m_i+1}}\right)^{[k:\Q]}
\end{align*}
(if $\G$ is $k$-split, then $s'=0$, $l = k$, $\D_{l/k} = 1$; in the non-split case we use the fact that $s'>2$).

Since for $i$ large enough $m_i!\gg (2\pi)^{m_i+1}$ it is clear that for large enough $r$
\begin{equation*}
\frac{1}{(\pi/12)(r+1)^3}\prod_{i=1}^{r}\frac{m_i!}{(2\pi)^{m_i+1}} > 1.
\end{equation*}
An easy direct computation shows that starting from $r=30$,
$$
10^{-2\epsilon'}\D_k^{\d - 2}\D_{l/k}
\left(\frac{1}{(\pi/12)(r+1)^3}\prod_{i=1}^{r}\frac{m_i!}{(2\pi)^{m_i+1}}\right)^{[k:\Q]} \ge
\D_k^{\d-2}\D_{l/k}.
$$
So for $r\ge 30$, $\delta = \d - 2$ and any field $k$ we have $B(\G/k) \ge \D_k^\delta \D_{l/k}$,
the finite set of the exceptional fields is empty and $(i)$ is proved.
\medskip

To proceed with the argument let us remark that
\begin{equation*}
B(\G/k) \ge \D_k^{\d-2}\D_{l/k}c,
\end{equation*}
where $c>0$ depends only on the absolute type of $\G$ and degree $d = [k:\Q]$. So, if the degree
$d$ is fixed, then for any $z>0$ which will be chosen later we have
\begin{equation}\label{325}
B(\G/k) \ge \D_k^{\d-2-z}\D_{l/k} \D_k^z c \ge \D_k^{\d-2-z}\D_{l/k},
\end{equation}
for all $k$ with $\D_k\ge c^{-z}$. Since there are only finitely many number fields with a bounded
discriminant, (\ref{325}) holds for all but finitely many $k$ of degree $d$. Since we always have
$\d > 2$, this allows to assume (at least when $\G$ is not $^6\Dn_4$) that the degree of $k$ is large enough.
\medskip

We now come to the case $(ii)$. Let $\G$ be still not of type $^6\Dn_4$. By the previous remark we can
suppose that $[k:\Q]$ is large enough. We have the following lower bound for $\D_k$ due
to A.~Odlyzko~\cite[Theorem~1]{Odlyzko}:
\begin{equation}\label{33}
{\rm If\ } [k:Q] > 10^5,{\rm\ then\ } \D_k \ge 55^{r_1(k)}21^{2r_2(k)}.
\end{equation}

So for $[k:Q] > 10^5$,
\begin{equation*}
B(\G/k) \ge 10^{-2\epsilon'}\left(\frac{21^{\d - 2 - \delta}}{(\pi/12)(r+1)^3}
\prod_{i=1}^{r}\frac{m_i!}{(2\pi)^{m_i+1}}\right)^{[k:\Q]} \D_k^\delta \D_{l/k}.
\end{equation*}
A direct case by case verification shows that for $\delta=1$ the latter expression is
$\ge \D_k^\delta \D_{l/k}$. So if we put $z = \d - 3$ in~(\ref{325}), then we obtain that
for all but finitely many $k$, $B(\G/k) \ge \D_k \D_{l/k}$.
\medskip

Let now $\G/k$ be a triality form of type $^6\Dn_4$. We have
\begin{gather*}
\epsilon = 2,\ \epsilon' = 1,\ n = 2,\ s' = 2.5,\ \{m_i\} = \{1,3,5,3\};\\
[l:k] = 3 {\rm\ and\ } a(l)\le 3[k:\Q].
\end{gather*}
So if $[k:Q] > 10^5$,
\begin{equation*}
B(\G/k) \ge 10^{-2}\left(\frac{21^{14-3 -\delta}}{(\pi/12)\cdot2^3}
\cdot\frac{6\cdot120\cdot6}{(2\pi)^{16}}\right)^{[k:\Q]} \D_k^\delta \D_{l/k}.
\end{equation*}
For $\delta = 1$ it is $\ge \D_k \D_{l/k}$. If $[k:\Q] \le 10^5$ we still have the inequality (\ref{325})
(the precise formula for the constant $c$ would be different but it is not essential), so for
all but finitely many $k$, again $B(\G/k)\ge \D_k \D_{l/k}$. The case $(ii)$ is now settled completely.
\medskip

Let $\G/k$ be of type $\An_3$ or $\Bn_2$. As before we can assume $[k:Q] > 10^5$. We have
$$
B(\G/k) \ge 10^{-2}\left(\frac{\pi}{12}\right)^{-[l:\Q]}
\left(\frac{21^{\d - [l:k] - \delta}}{n^{1+[l:k]}}
\prod_{i=1}^{r}\frac{m_i!}{(2\pi)^{m_i+1}}\right)^{[k:\Q]} \D_k^\delta \D_{l/k}.
$$
Now, $n = 4$ and $n = 2$ for the types $\An_3$ and $\Bn_2$, respectively; if $l\neq k$, then
$[l:\Q] = 2[k:\Q]$. Using this it is easy to check that if $\delta = 0.1$,
then $B(\G/k)\ge \D_k^\delta \D_{l/k}$ in each of the possible cases.
\medskip

It remains to consider $(iv)$. This is the most difficult case but the proof almost
repeats the argument of~\cite[Prop.~6.1(vi)]{BP}.

With the notations of~\cite{BP}, for $[l:\Q] > 10^5$ we have
\begin{align*}
\lefteqn{
B(\G/k) = \D_k^4 \cdot 3^{- a(k)- a(l)}h_l^{-1}
(\D_l/\D_k^{[l:k]})^{5/2}\left(\frac{1}{2^4\pi^5}\right)^{[k:\Q]}} \\
 &\quad = \D_k^{\delta} \D_l^{2-\delta/2} (\D_{l/k})^{1/2} 3^{- a(k)- a(l)}h_l^{-1}
          \left(\frac{1}{2^4\pi^5}\right)^{[k:\Q]} \\
 &\quad \ge \D_k^{\delta}(\D_{l/k})^{1/2} \frac{0.02}{s(s-1)}
 \left(\frac{55^{(4-\delta-s)/2}}{2\cdot3^{3/2}\cdot\pi^{(6-s)/2}}\right)^{r_1(l)}
            \left(\frac{21^{(4-\delta-s)}}{2^{(4-s)}\cdot3^{2}\cdot\pi^{(5-s)}}\right)^{r_2(l)} \\
 &\quad\quad \cdot\exp\big((3-\delta-s)Z_l(s) - (4-\delta-s)(\frac{c_1}{2} + (s-1)^{-1})\\
 &\phantom{xxxxxxxxxxxxxxxxxxxx}+ (0.1-(c_3 + c_4)(s-1))a(l)\Big).
\end{align*}
Now let $\delta = 0.01$. Since
\begin{align*}
 & 55^{2.99/2}(2\cdot3^{3/2}\cdot\pi^{5/2})^{-1} > 2.19,\quad 21^{2.99}(2^{3}\cdot3^{2}\cdot\pi^{4})^{-1} > 1.28,\\
 {\rm and}\ & \exp((3-\delta-s)Z_l(s))\ge 1\ {\rm if}\ s<2-\delta,
\end{align*}
by choosing $s>1$ sufficiently close to $1$, we obtain that there is an absolute constant $c_6$
such that
$$
\D_k^4 \cdot 3^{- a(k)- a(l)}h_l^{-1} (\D_l/\D_k^{[l:k]})^{5/2}\left(\frac{1}{2^4\pi^5}\right)^{[k:\Q]}
%&\quad\ge
\ge \D_k^{0.01}(\D_{l/k})^{1/2} 2.19^{r_1(l)} 1.28^{r_2(l)} c_6.
$$
The right-hand side is $\ge \D_k^{0.01}(\D_{l/k})^{1/2}$ if $[l:Q]$ is large enough, say, $[l:Q] > d_l$ (and $d_l \ge 10^5$).

\noindent
If $[l:\Q] < d_l$, then $[k:\Q] < d_l$ and by $(3)$ for all but finitely many fields $k$
we have
\begin{align*}
B(\G/k) \ge \D_k \D_{l/k} \ge \D_k^{0.01} \D_{l/k}^{0.5}.
\end{align*}
\end{proof}

\begin{remark}
The proof provides explicit values of $\delta_1$, $\delta_2$ for each of the types,
however in many cases the bound for $B(\G/k)$ can be improved. This will require more
careful argument and can be useful for particular applications.
\end{remark}

\section{Proof of the theorem. The upper bound}\setcounter{equation}{0}

As before $H$ denotes a connected semi-simple Lie group whose almost simple factors
are all non-compact and have the same type different from $\An_1$, $\G$ is an absolutely almost simple
simply connected $k$-group admissible in the sense that there exists a continuous surjective
homomorphism \mbox{$\G(k\otimes_\Q\R)^o\to H$} with a compact kernel.

\subsection{Counting number fields}\label{s41}
For a (maximal) arithmetic subgroup $\Gamma$ of $H$ we have (see~\ref{2:BP})
\begin{equation}\label{41}
\mu(H/\Gamma) \ge  \frac12 B(\G/k)\tau_k(\G)\F
\end{equation}
where \quad $k$ is the field of definition of $\Gamma$,
\begin{gather*}
\intertext{\quad\quad\quad\ $\G/k$ is a $k$-form from which $\Gamma$ is induced (see~\ref{s21});}
B(\G/k) = \D_k^{\d} n^{-\epsilon a(k)-\epsilon' a(l)}h_l^{-\epsilon'}
(\D_l/\D_k^{[l:k]})^{s'}\left(\prod_{i=1}^{r}\frac{m_i!}{(2\pi)^{m_i+1}}\right)^{[k:\Q]}
\!\!\!\!\!\!\!\!\!;\\
\tau_k(\G) = 1;\\
\F = \prod_{v\in V_f} f_v > 1 {\rm\ \ will\ be\ considered\ later}.\\
\end{gather*}
By Proposition~\ref{prop:BP} for all but finitely many number fields $k$
\begin{equation*}\mu(H/\Gamma) \ge c_1 \D_k^{\delta_1} \D_{l/k}^{\delta_2},\end{equation*}
where $\delta_1$, $\delta_2$ are the constants determined by the absolute type of $\G$ (which is the
type of almost simple factors of $H$).

So for large enough $x$, if $\mu_S(H/\Gamma) < x$, then $\D_k < (x/c_1)^{1/\delta_1}$,
$\D_{l/k} < (x/c_1)^{1/\delta_2}$.

By Proposition~\ref{prop:ND}(ii), the number of such fields $k$ is at most
\begin{equation*}(x/c_1)^{\beta((x/c_1)^{1/{\delta_1}})}\le x^{c_2\beta(x)}, \end{equation*}
and by Proposition~\ref{prop:ND}(i), for each $k$ the number of such extensions $l$ is at most
\begin{equation*}
c\D_k^{b_1}(x/c_1)^{b_2/{\delta_2}} \le c(x/c_1)^{b_1/\delta_1}(x/c_1)^{b_2/\delta_2}\le x^{c_3}.
\end{equation*}
It follows that the number of all admissible pairs $(k,l)$ is bounded by
\begin{equation*}x^{c_2\beta(x)+c_3},\end{equation*}
and, moreover, since $k\neq\Q$ implies $[k:\Q]\le c\log\D_k$, for all admissible $k$ we have
\begin{equation*} a(k) \le c_4\log x.\end{equation*}

\subsection{Non-cocompact case}\label{s42}
If $\Gamma$ is non-cocompact, the degree of the field of definition of $\Gamma$ is bounded. Indeed,
the non-cocompactness of $\Gamma$ implies that the corresponding algebraic group $\G$ is $k$-isotropic,
so $\G/k_v$ is non-compact for every $v\in V$. It follows that the number of infinite places of $k$ is
equal to the number $\# S$ of almost simple factors of $H$, so $[k:\Q]\le 2\# S$.

Now in~\ref{s41} we can consider only the number fields $k$ with $[k:\Q]\le 2\# S$ and the
number fields $l$ with $[l:\Q]\le 3[k:\Q]\le 6\# S$. By Proposition~\ref{prop:ND}(i), for large enough
$x$ the number of admissible pairs $(k,l)$ is at most
\begin{equation*} x^{c_5},\ c_5 = c_5(\# S) \end{equation*}
(in fact, here we could use a weaker result by W.~M.~Schmidt,
{\em Ast\'erisque}, {\bf 228}(1995), 189--195, who showed that the number
of degree $n$ extensions $l$ of $k$ with $\D_{l/k} < x$ is bounded by $C(n)x^{(n+2)/4}$).

\subsection{Counting $k$-forms}\label{s43}
Given an admissible pair $(k,l)$ of number fields, there exists
a unique quasi-split $k$-form $\Gqs$ for which $l$ is the splitting field (or a certain subfield of
the splitting field if $\Gqs$ is of type $^6\Dn_4$ and $[l:\Q] = 3$). So we have an upper bound
for the number of quasi-split groups for which there can exist an inner form that defines an
arithmetic subgroup of covolume less than $x$. We now fix a quasi-split $k$-form $\Gqs$ and estimate
the number of admissible inner forms. Since every inner equivalence class of $k$-forms contains
a unique quasi-split form, this will give us a bound on the total number of admissible $\G/k$.
\medskip

By the assumption, $\prod_{v\in V_\infty(k)} \G(k_v)$ is isogenous to $H\times K$ ($K$ is a compact
Lie group), so the $k_v$-form of $\G$ is almost fixed at the infinite places of $k$. More precisely,
let $c_h$ be the number of non-isomorphic almost simple factors of $H$. For each $v\in V_\infty(k)$,
$\G(k_v)$ is isomorphic to one of $c_h$ non-compact (simply connected) groups or is compact, and the
number of places $v$ at which $\G(k_v)$ is non-compact $n_h = \# S$. This implies that the number of
variants for $\G(k_v)$ at the infinite places of $k$ is bounded by
\begin{equation*}
c_h^{n_h} \binom{a(k)}{n_h} < (c_h a(k))^{n_h} \le (\log x)^{c_6}
\end{equation*}
($\binom{\cdot}{\cdot}$ denotes the binomial coefficient).
\medskip

Let now $v$ be a finite place of $k$. The inner $k_v$-forms of $\G$ correspond to the elements of the
first Galois cohomology set $\H^1(k_v, {\overline \G})$, ${\overline \G}$ is the adjoint group of $\G$.
The order of $\H^1(k_v, {\overline \G})$ can be computed from the cohomological exact sequence
\begin{equation*}
\H^1(k_v, {\G}) \to \H^1(k_v, {\overline \G}) \stackrel{\delta}{\to} \H^2(k_v, {\rm C}),
\end{equation*}
which corresponds to the universal $k_v$-covering sequence of groups
\begin{equation*}
1 \to {\rm C} \to {\G} \to {\overline \G} \to 1.
\end{equation*}
For a simply connected $k_v$-group ${\G}$ the first cohomology
$\H^1(k_v, {\G})$ is trivial by a theorem of Kneser~\cite{Kn}, so
$\delta$ is injective. Furthermore, the group $\H^2(k_v, {\rm C})$ can be identified
with a subgroup of the Brauer group of $k_v$ and then explicitly computed using
results from the local class field theory. We refer to~\cite[Chapter~6]{PlR}
for details and explanations. As a corollary here we have that the number of inner
$k_v$-forms is bounded by $n^\epsilon$ in the notation of~\ref{1:BP} (recall
that $n^{\epsilon} = \# {\rm C}$ is the order of the center of $\G$).

Let $T_1\subset V_f(k)$ be a (finite) subset of the nonarchimedean places of $k$ such that $\G$ is not
quasi-split over $k_v$ for $v\in T_1$. It follows from~\cite[Prop.~2.10]{P} that there exists a constant
$\delta > 0$, which depends only on the absolute type of $\G$, such that for every $v\in T_1$,
\begin{equation}\label{43}
f_v \ge n^{-\epsilon}(\#\Xi_{\Theta_v})^{-1}e_v \ge q_v^\delta
\end{equation}
($q_v$ denotes the order of the residue field of $k$ at $v$).

Indeed, we can take $\delta = \log(2^{r_v}n^{-\epsilon})$ if the absolute type of $\G$ is not $\An_2$
and $\delta = \log(2^2\cdot3^{-1}) = 0.415\dots$ for the type $\An_2$, and then check that $\delta>0$ and
inequality (\ref{43}) holds going through the case-by-case consideration in~\cite[Appendix C.2]{BP}.

To a set $T\subset V_f(k)$ we can assign an ideal $\I_T$ of $\cO_k$ given by the product of prime ideals
defining the places in $T$. Reversely, each squarefree ideal of $\cO_k$ uniquely defines a subset
$T$ in $V_f(k)$ corresponding to its prime decomposition. Note also that $\prod_{v\in T}q_v = {\rm Norm}(\I_T)$.

Now for an arithmetic subgroup $\Gamma$ induced from $\G$ we have
\begin{equation*}
\mu(H/\Gamma) = \frac12 B(\G/k)\tau_k(\G)\F \ge c_7\prod_{v\in T_1} q_v^\delta.
\end{equation*}
This implies that if $\mu(H/\Gamma) \le x$, then ${\rm Norm}(\I_{T_1}) = \prod_{v\in T_1}q_v \le x^{c_8}$.
By Proposition~\ref{prop:NP}(ii) the number of variants for $T_1$ is bounded by $x^{c_9}$,
moreover, since for every $v\in V_f$, $q_v\ge 2$, for every such a set $T_1$ we have
$\# T_1 \le c_{10}\log x$.
\medskip

Now the Hasse principle implies that a $k$-form of $\G$ is uniquely determined
by $(\G(k_v))_{v\in V(k)}$. The Hasse principle for semi-simple groups is valid
due to the work of Kneser, Harder, Chernousov (cf.~\cite[Chapter~6]{PlR}). So
the number of the admissible $k$-forms is at most
\begin{equation*}
(\log x)^{c_6} x^{c_8} n^{\epsilon c_{10}\log x} \le x^{c_{11}}.
\end{equation*}

\subsection{Counting collections of parahorics} \label{s44}
For a given large enough $x$ we have defined a
collection of $\G/k$ for which there exists a (centrally) $k$-isogenous
group $\G'$ which may give rise to the arithmetic subgroups $\Gamma\subset H$
with $\mu(H/\Gamma) < x$. The number of such $k$-groups $\G$
is finite and can be bounded as in~\ref{s43}, but still each $\G/k$ defines countably
many maximal arithmetic subgroups. We shall now fix a group $\G/k$ and estimate the
number of coherent collections of parahoric subgroups of $\G$ which can give rise to
the maximal arithmetic subgroups with covolumes less than $x$. In the classical
language, what we are going to do in this section is to count the number of
admissible genera.

We use again the local to global approach. Let us fix a central $k$-isogeny $i: \G\to \G'$ with
$\G'$ such that $\G'_S$ projects onto $H$. Every
maximal arithmetic subgroup $\Gamma\subset \G'_S$ is associated to some coherent collection
$P = (P_v)_{v\in V_f}$ of parahoric subgroups of $\G$ (see~\cite[Prop.~1.4]{BP}):
\begin{equation}
\Gamma = N_{\G'}(i(\Lambda)),\ \Lambda = \G(k)\ \cap \!\!\prod_{v\in V_f(k)} P_v.
\end{equation}
The image of $\Gamma$ in $H$ is an arithmetic subgroup and every maximal arithmetic subgroup
of $H$ can be obtained as a projection of some such $\Gamma$.

For almost all finite places $v$ of $k$, $\G$ is quasi-split over $k_v$ and splits over
an unramified extension of $k_v$. Moreover, for almost all such $v$, $P_v$ is hyperspecial.
Any two hyperspecial parahoric subgroups of $\G(k_v)$ are conjugate under the action
of the adjoint group $\Gb(k_v)$ \cite[Sect. 2.5]{T}, so $P$ is determined
up to the action of $\Gb(\A_f)$ by the types of $P_v$ at the remaining places. Using this
we shall now count the number of $P$'s.

Let as in~\ref{s43}, $T_1$ denote the set of places of $k$ for which $\G$ is not quasi-split.
By the previous argument we have
\begin{equation*} \# T_1 \le c_{10}\log x,\ \
\# ({\rm variants\ for}\ T_1) \le x^{c_9}
\ \ \rm{(see\ \ref{s43})}.\end{equation*}
Let $R$ denote the set of places for which $\G$ is quasi-split but is not split over
an unramified extension of $k_v$. For such places $v\in V_f$, $l_v = l\otimes_k k_v$ is a
ramified extension of $k_v$ and so by the formula from \cite[Appendix]{P}, each of
such places contributes to $\D_{l/k}$ a power of $q_v$. Using again
Proposition~\ref{prop:NP}(ii) and the inequality $\D_{l/k}\le x^c$ from~\ref{s41},
we obtain
\begin{equation*}
\# R \le c_{12}\log x,\ \
\# ({\rm variants\ for}\ R) \le x^{c_{13}}.
\end{equation*}
Finally, let $T_2\subset V_f\backslash(T_1\cup R)$ be the set of places for which $P_v$ is
not hyperspecial. If $v\in T_2$, then by~\cite[Prop.~2.10(iv)]{P}
\begin{equation*} e_v \ge (q_v+1)^{-1}q_v^{r_v+1}.\end{equation*}
Similarly to~(\ref{43}) it implies
\begin{equation*} f_v \ge q_v^\delta.\end{equation*}
By Proposition~\ref{prop:NP}(ii) and the volume formula,
\begin{equation*}
\# T_2 \le c_{14}\log x,\ \
\# ({\rm variants\ for}\ T_2) \le x^{c_{15}}.
\end{equation*}

Now for a given $v\in V_f$ the number of the possible types of parahoric subgroups (parametrized
by the subsets of the set of simple roots) is bounded by a constant $c_t$ which depends
only on the absolute type of $\G$. We conclude that for a given $\G$ the number of $P$'s,
up to the action of $\Gb(\A_f)$, is at most
$$
c_t^{\#(T_1 \cup R \cup T_2)}\# ({\rm variants\ for}\ T_1 \cup R \cup T_2)
\le c_t^{(c_{10}+c_{12}+c_{14})\log x} x^{c_9+c_{13}+c_{15}} = x^{c_{16}}.
$$

\subsection{Counting conjugacy classes} \label{s45}
In this final step we give an upper bound for the number of conjugacy classes of arithmetic
subgroups associated to a fixed group $\G'/k$ and a given $\Gb(\A_f)$-orbit of collections
of parahoric subgroups $P$ of a simply connected group $\G$ centrally $k$-isogenous to $\G'$.
We are interested in the $\Gb(k)$-conjugacy classes of maximal subgroups associated to $P$
%which in turn amounts to conjugacy under the action of the adjoint group $\Gb(k)$. The $\Gb(k)$-conjugacy classes corresponding to $P$
which are indexed by the double cosets $\Gb(k)\backslash\Gb(\A)/\Gb_\infty\Pb$, where
$\Gb_\infty = \prod_{v\in V_\infty} \Gb(k_v)$, $\Pb_v$ is the stabilizer of $P_v$ in $\Gb(k_v)$
and $\Pb = \prod_{v\in V_f} \Pb_v$ is a compact open subgroup of
$\Gb_f = \prod_{v\in V_f}\Gb(k_v)$ (see \cite[Prop.~3.10]{BP}).
The number $c(\Pb)$ of the double cosets is called the {\it class number} of $\Gb$ with
respect to $\Pb$. The argument is similar to the proof of
Proposition~3.9 in \cite{BP} except that we need to get an explicit upper bound for
$c(\Pb)$.

Let $\omega$ be a non-zero invariant exterior $k$-form of top degree on $\Gb$;
such a form is unique up to multiplication by an element of $k^*$ and is called a
Tamagawa form. We denote by $|\omega|$ the Haar measure on the ad\`ele group $\Gb(\A)$
determined by $\omega$. % and also corresponding Haar measures on $\Gb(k_v)$, $v\in V(k)$.
The natural embedding of $k$ into $\A$ gives an embedding of $\Gb(k)$ in $\Gb(\A)$ and
it is well known that the image of $\Gb(k)$ is a lattice in $\Gb(\A)$. By the product
formula its covolume with respect to the measure $|\omega|$ does not depend on the
choice of the form $\omega$, thus the number
$\tau_k(\Gb) := \D_k^{-\dim(\Gb)/2} |\omega|(\Gb(k)\backslash\Gb(\A))$ is correctly
defined. It is called the {\it Tamagawa number} of $\Gb/k$.
By a theorem of T.~Ono \cite{Ono}, $\tau_k(\Gb)$ is bounded by a constant multiple of the
order of the center of the simply connected covering group $\G$ multiplied by $\tau_k(\G)$. The
Tamagawa number of a simply connected group is equal to $1$ according to the Weil conjecture
which has been proved completely for the groups over number fields due to the work of many
people (see \cite[3.3]{P} for a short discussion). Therefore we have
\begin{equation}\label{eq_tau}
|\omega|(\Gb(k)\backslash\Gb(\A)) = \tau_k(\Gb) \D_k^{\dim(\Gb)/2}\le c_{17} \D_k^{\dim(\Gb)/2}
\end{equation}
where $c_{17}$ depends only on the absolute type of $\G$.

Coming back to the problem of bounding the class number $c(\Pb)$, we recall that the
double cosets $\Gb(k)\backslash\Gb(\A)/\Gb_\infty\Pb$ correspond bijectively to the
orbits of $\Gb_\infty\Pb$ on $\Gb(k)\backslash\Gb(\A)$ which are open. Given an
upper bound for $|\omega|(\Gb(k)\backslash\Gb(\A))$, in order to give a bound for
$c(\Pb)$ it is enough to obtain a uniform lower bound for the $|\omega|$-volumes of these
orbits. The double cosets are represented by elements of $\Gb_f$, so it is sufficient
to consider the orbit of the image of $a\in\Gb_f$ which is isomorphic to
$\Gamma_a\backslash\Gb_\infty a\Pb a^{-1}$, $\Gamma_a = \Gb(k)\cap\Gb_\infty a\Pb a^{-1}$.
Let $\Gamma_a'$ be the projection of $\Gamma_a$ to $\Gb_\infty$ with respect to the
decomposition $\Gb(\A) = \Gb_\infty\times\Gb_f$.  As $a\Pb a^{-1}$
is a compact open subgroup of $\Gb_f$, $\Gamma_a'$ is an arithmetic subgroup
of $\Gb_\infty$. We have
\begin{equation*}
|\omega|(\Gamma_a\backslash\Gb_\infty a\Pb a^{-1}) =
|\omega|_\infty(\Gamma_a'\backslash\Gb_\infty) |\omega|_f(\Pb),
\end{equation*}
where $|\omega|_\infty$, $|\omega|_f$ denote the product measures on
$\Gb_\infty$, $\Gb_f$ corresponding to $\omega$.

In order to estimate the factors in the right-hand side of the formula, for each
$v\in V(k)$ we shall relate the measure $|\omega|$ to the canonical measure
$|\omega_{\Gb_v}|$ on $\Gb(k_v)$ defined in \cite[Sections~4, 11]{Gross}. In
particular, if $\G$ is simply connected the measure $|\omega_{\G_v}|$ coincides
with the measure $\gamma_v\omega_v^*$ which is used for the local computations in
\cite{P}; for $v\in V_\infty(k)$, $|\omega_{\Gb_v}|$ is equal to the measure
$\mu$ on $\Gb(k_v)$ defined as in \ref{s25}; and for all but finitely many $v$,
$|\omega_{\Gb_v}| = |\omega|_v$. Let $\gamma_v$ denote the ratio
$|\omega_{\Gb_v}|/|\omega|_v$ which by the previous remark is equal to $1$ for
all but finitely many places $v$. Hence
\begin{equation*}
|\omega|(\Gamma_a\backslash\Gb_\infty a\Pb a^{-1}) =
\mu_\infty(\Gamma_a'\backslash\Gb_\infty) \prod_{v\in V_f} |\omega_{\Gb_v}|(\Pb_v) / \prod_{v\in V}\gamma_v.
\end{equation*}

We now recall the main result of \cite{BP}, which implies that covolumes of arithmetic subgroups of
$\Gb_\infty$ with respect to the measure $\mu$ are bounded from below by a universal constant and thus
$\mu_\infty(\Gamma_a'\backslash\Gb_\infty) \ge \mu_0$.

The crucial ingredient which allows us to carry out the required estimates is the product formula
for $\gamma_v$. It was obtained in \cite[Theorem~1.6]{P} for the simply connected groups and later
extended by B.~Gross to arbitrary reductive groups defined over number fields (see also \cite{Ku}
for the groups over global function fields). Thus by \cite[Theorem~11.5]{Gross} we have
\begin{equation*}
\prod_{v\in V} \gamma_v = (\D_l/\D_k^{[l:k]})^{\frac12s}
\left(\prod_{i=1}^{r}\frac{m_i!}{(2\pi)^{m_i+1}}\right)^{[k:\Q]}\!\!\!\!\!\!\!\!\!.
\end{equation*}

Finally, we shall make use of the following inequality:
\medskip

{\it Claim.} $|\omega_{\Gb_v}|(\Pb_v) \ge |\omega_{\G_v}|(P_v) = e(P_v)^{-1}$.
\medskip

The proof of this claim which is given below is quite technical but not conceptually new, related
questions were studied in detail and full generality in \cite{Gross}, \cite{Ku}. The argument
falls into several steps.

\def\Gx{\underline{\mathrm G}_x}
\def\Gix{{\underline{\mathrm G}'}_x^0}
\def\Gxs{\overline{\mathrm G}_x}
\def\Gixs{{\overline{\mathrm G}'}_x^0}
\def\og{|\omega_\G|}
\def\ogi{|\omega_{\G'}|}
\def\F{{\mathrm F}}
\def\Bb{\overline{B}}
\def\Go{\underline{\mathrm G}_\Omega}
\def\Gio{{\underline{\mathrm G}'}_\Omega^0}
\def\Gos{\overline{\mathrm G}_\Omega}
\def\Gios{{\overline{\mathrm G}'}_\Omega^0}

Let $K = k_v$ be a nonarchimedean local field, $\cO$ its ring of integers, $\G$ a
simply connected semi-simple $K$-group, $i: \G \to \G'$ a central $K$-isogeny (we actually need only
the case $\G'=\Gb$) and let $X = X(\G)$ denote the Bruhat-Tits building of $\G/K$.

{\bf 1.}
We shall assume first that the groups $\G$ and $\G'$ are quasi-split over $K$. Let $x\in X$ be a special
vertex in $X$ chosen as in [Gr, Sect.~4] (see also [P, 1.2]). The Bruhat-Tits theory assigns
to $\G'/K$ and $x\in X(\G)$ a smooth affine group scheme $\Gix$ over $\cO$. Its generic fiber
is isomorphic to $\G'/K$ and its special fiber $\Gixs$ is connected. Let $P_x = \Gx(\cO)
(= \Gx^0(\cO))$, $P'_x = \Gix(\cO)$. Then $P_x$ (resp. $P'_x$) is an open compact subgroup
of $\G(K)$ (resp. $\G'(K)$), $P_x$ is the stabilizer of $x$ in $\G(K)$ and $P'_x$ is contained
in the stabilizer of $x$ in $\G'(K)$ with finite index. Recall also that the measure
$\og$ (resp. $\ogi$) corresponds to a differential $\omega_\G$ (resp.
$\omega_{\G'}$) of top degree on $\G$ (resp. $\G'$) over $K$ which has good
reduction (see [Gr, Sect.~4]). This brings us to the conditions of
Proposition~I.2.5 of [Oe], which implies
\begin{equation*}
\og(P_x) = \#\Gxs(\F_q) q^{-\mathrm{dim}\:\G},\
\ogi(P'_x) = \#\Gixs(\F_q) q^{-\mathrm{dim}\:\G'},
\end{equation*}
($\F_q$ denotes the residue field of $K$).

Since $\G$ and $\G'$ are isogenous, $\Gxs$ and $\overline{\mathrm G}'_x$ are
isogenous. Hence, $\mathrm{dim}\:\G = \mathrm{dim}\:\G'$,
and by Lang's theorem, $\#\Gxs(\F_q) = \#\Gixs(\F_q)$.  Thus we obtain
$\og(P_x) = \ogi(P'_x)$.

{\bf 2.}
Let now $C$ be a chamber of $X$ which contains $x$ and let $\Omega$ be a subset
of $C$. Denote by $I_C$ (resp. $I'_C$) the Iwahori subgroup of $\G(K)$ (resp.
$\G'(K)$) corresponding to $C$. Note that by definition $I'_C$ is the preimage
in $\G'(K)$ of a Borel subgroup $\Bb'$ of $\Gixs(\F_q)$. Let $P_\Omega$ (resp.
$P'_\Omega$) be the parahoric subgroup of $\G(K)$ (resp. $\G'(K)$) associated to
$\Omega$; so $P_\Omega = \Go(\cO)$, $P'_\Omega = \Gio(\cO)$ and any parahoric
subgroup of $\G(K)$ is conjugate to some $P_\Omega$. The inclusion
$\Omega\subset C$ induces a group scheme homomorphism $\rho_{\Omega C}:
{\underline{\mathrm G}'}^0_C \to \Gio$ whose reduction maps the group
${\overline{\mathrm G}'}^0_C$ onto a Borel subgroup $\Bb'$ of $\Gios$.
Therefore we have
\begin{equation*}
[P'_\Omega:I'_C] = [\Gios(\F_q):\Bb'] = [\Gos(\F_q):\Bb] = [P_\Omega:I_C],
\end{equation*}
as $\Gios$ is isogenous to $\Gos$, $\Bb'$ is isogenous to $\Bb$ and all the
groups are connected. It follows that
$|\omega_{\G'}|(P'_\Omega) = |\omega_{\G}|(P_\Omega)$.

{\bf 3.}
We finally note that $P'_\Omega \subset \overline{P}'_\Omega$
($\overline{P}'_\Omega$ denotes the stabilizer of $\Omega$ in $\G'(K)$)
and thus
$|\omega_{\G'}|(\overline{P}'_\Omega) \ge |\omega_{\G'}|(P'_\Omega) = |\omega_{\G}|(P_\Omega)$,
which implies the desired inequality in the quasi-split case.

{\bf 4.}
In order to extend this result to the general case we have to recall the definition of
the canonical measure $\og$ for the general $\G$ by pull-back from the quasi-split inner
form (see \cite[p.~294]{Gross}) and its interpretation in terms of the volume form $\nu_G$
associated to an Iwahori subgroup of $\G(K)$ described in \cite[pp.~294--295]{Gross}. The
latter allows us to apply the argument similar to step 1 to Iwahori subgroups $I_C$ and $I'_C$
corresponding to a chamber $C$ of $X(\G)$,
proving $\og(I_C) = \ogi(I'_C)$. All the rest of the proof does not depend on the quasi-split
assumption and the claim follows.

\medskip

%Indeed, first of all there is
%a standard way to reduce the question to the case when $\G$ and $\Gb$ are both quasi-split over
%$k_v$~--- see e.g. the proof of Prop.~4.7 in \cite{Gross}. Now by the definition of measures $|\omega_{\G_v}|$
%and $|\omega_{\Gb_v}|$, Prop.~I.2.5 of \cite{Oe} and Lang's theorem, which says that connected isogenous
%algebraic groups defined over a finite field ${\mathrm F}$ have the same number of ${\mathrm F}$-rational
%points, we obtain $|\omega_{\G_v}|(I_v) = |\omega_{\Gb_v}|(I'_v)$, where $I_v$ is an Iwahori subgroup
%of $\G(k_v)$ and $I'_v$ is the corresponding Iwahori subgroup of the (non simply connected) group
%$\Gb(k_v)$ defined as in \cite[3.7]{T}. Note that $I_v'$ is contained in $\overline{I}_v$ with a finite index
%and, moreover, $[P_v:I_v] = [P'_v:I'_v]$, where $P'_v$ is a maximal subgroup of $\Pb_v$ whose reduction
%is contained in the neutral component of the special fiber of the Bruhat--Tits group scheme associated
%with $\Pb_v$. In particular, $P'_v \subset \Pb_v$ and the desired inequality follows. Related questions
%are considered in detail and full generality in \cite{Gross}, \cite{Ku}.

Let us collect together the results of this section. We obtain:
\begin{align}\label{eq_cl}
\lefteqn{
c(\Pb) \le \frac{|\omega|(\Gb(k)\backslash\Gb(\A)) \prod_{v\in V}\gamma_v}
{\mu_0 \prod_{v\in V_f}e(P_v)^{-1}} } \notag\\
&\quad \le \frac{1}{\mu_0}\D_k^{\dim(\Gb)/2}(\D_l/\D_k^{[l:k]})^{\frac12s}
\left(\prod_{i=1}^{r}\frac{m_i!}{(2\pi)^{m_i+1}}\right)^{[k:\Q]}\tau_k(\Gb) \E(P).
\end{align}
This formula can be viewed as an extension of the upper bound for the class number from \cite[Theorem~4.3]{P}.

We now bound the right-hand side of (\ref{eq_cl}).
By \ref{s41}, \ref{s42} we have $\D_k < (x/c_1)^{1/\delta_1}$,  $\D_{l/k} < (x/c_1)^{1/\delta_2}$ and
$[k:\Q]\le c\log\D_k$. By (\ref{eq_tau}), $\tau_k(\Gb)\le c_{17}$. From \ref{s43} and \ref{s44} it follows that
if $\mu(H/\Gamma) \le x$, then $\prod_{v\in V_f} e_v \le x^{c_{18}}$ for some constant $c_{18}$ which depends
only on the type of almost simple factors of $H$. Hence it follows from (\ref{eq_cl}) that there exists a
constant $c_{19}$ such that
$$ c(\Pb) \le  x^{c_{19}}.$$

\subsection{The upper bounds}
It remains to combine the results of the previous sections to get the upper bounds.
By \ref{s41}, \ref{s43}, \ref{s44}, \ref{s45}
\begin{equation*}
m^u_H(x) \le x^{c_2\beta(x)+c_3} x^{c_{11}} x^{c_{16}} x^{c_{19}} \le x^{B\beta(x)},
\end{equation*}
and constant $B$ depends only on the type of almost simple factors of $H$.
\medskip

By \ref{s42}, \ref{s43}, \ref{s44}, \ref{s45}
\begin{equation*}
m^{nu}_H(x) \le x^{c_5} x^{c_{11}} x^{c_{16}} x^{c_{19}} \le x^{B'},
\end{equation*}
constant $B'$ depends on the type and the number of almost simple factors of $H$.

\section{Proof of the theorem. The lower bound}\setcounter{equation}{0}

\subsection{Cocompact case.}\label{s51}
A theorem of Borel and Harder~\cite{BH} implies that a semi-simple group over a local field
of characteristic $0$ contains cocompact arithmetic lattices. The method of~\cite{BH} actually
proves the existence of such lattices defined over a given filed $k$, which satisfies a
natural admissibility condition, for any isotypic semi-simple Lie group. So if $H$ has $a_1$
real and $a_2$ complex almost simple factors (all of the same type) and $k$ is a number
field with $>a_1$ real and precisely $a_2$ complex places, then $H$ contains a cocompact
arithmetic subgroup $\Gamma_1$ defined over $k$.

Let $\Gamma_0$ be a maximal arithmetic subgroup of $H$ which contains $\Gamma_1$. There exists
an absolutely almost simple simply connected $k$-group $\G$ and  a principal arithmetic subgroup
$\Lambda_0$ of $\G$ such that $\Gamma_0 = N_H(\phi(\Lambda_0))$.

We assume $x$ is large enough and estimate the number of principal arithmetic
subgroups $\Lambda\subset \G(k)$ which are associated to the coherent collections of parahoric
subgroups of $\cO$-maximal types (\cite{R}, \cite{CR}) and such that $\mu_S(\G_S/\Lambda) < x$.
Then for $\Gamma = N_H(\phi(\Lambda))$ we also have $\mu_S(\G_S/\Gamma) < x$. Moreover,
by Rohlfs' theorem each such $\Gamma$ is a maximal arithmetic subgroup of $H$ and all maximal
arithmetic subgroups of $H$ are obtained as the normalizers of the images of the principal
arithmetic subgroups corresponding to $\cO$-maximal collections of parahorics.

The condition of maximality for the type of a collection of parahoric subgroups
$P = (P_v)_{v\in V_f}$ is a local condition on the types of $P_v$ at each $v\in V_f$,
while $\cO$-maximality requires an additional global restriction which is
needed to further narrow down the set of admissible collections of parahoric subgroups of maximal
types. We shall not give precise definitions here referring the reader to the above cited
papers. What is important for our argument is that given $P_0 = (P_{0,v})_{v\in V_f}$, a
collection of parahoric subgroups of $\cO$-maximal type, for every $v_0\in V_f$ there exists another
$\cO$-maximal collection $P = (P_v)_{v\in V_f}$ such that for $v\neq v_0$, $P_v = P_{0,v}$
and $P_{v_0}\not\cong P_{0,v_0}$. This is clearly true: for the groups of the absolute rank greater
than one (which is our standing assumption) it is enough to consider the maximal types corresponding
to single vertices of the affine Dynkin diagram and for such types $\cO$-maximality can be easily
checked.
%\medskip

We have
\begin{align*}
\mu_S(\G_S/\Lambda) & = \D_k^{\d}(\D_l/\D_k^{[l:k]})^{\frac12s}
\left(\prod_{i=1}^{r}\frac{m_i!}{(2\pi)^{m_i+1}}\right)^{[k:\Q]} \tau_k(\G)\:\E(P) \\
& = c_1\prod_{v\in T} e(P_v)/e(P_{0,v}) \\
& \le c_1 \prod_{v\in T} e(P_v), \ \ c_1 = \mu_S(\G_S/\Lambda_0),
\end{align*}
where $P_v$ (resp. $P_{0,v}$) is the closure of $\Lambda$ (resp. $\Lambda_0$) in $\G(k_v)$, $v\in V_f$; $T$
is a finite subset of the nonarchimedean places of $k$ for which $P_v\not\cong P_{0,v}$;
the constant $c_1$ depends on $\G/k$ and $\Lambda_0$ but does not depend on the choice of
$\Lambda$.
% The dependence of $c_1$ on $\Lambda_0$ can be further eliminated by looking at the construction
% of $\Lambda_0$ in \cite{BH}.

If $\prod_{v\in T} e(P_v) < x/c_1$, then $\mu_S(\G_S/\Lambda) < x$. There
exists a constant $\delta$ determined by the absolute type of $\G$ such that
for every $v\in V_f$ and every  parahoric subgroup $P_v\subset\G(k_v)$, $e(P_v)
\le q_v^{\delta}$ (e.g., take $\delta = {\rm dim}(\G)$). This implies
\begin{equation*}
\prod_{v\in T} e(P_v) \le \prod_{v\in T} q_v^{\delta}.
\end{equation*}
Hence $\prod_{v\in T} q_v < (x/c_1)^{1/\delta}$ is sufficient for $\mu_S(\G/\Lambda) < x$.
The number of variants for such sets $T$ is controlled via Proposition~\ref{prop:NP}(i)
(note that the field $k$ is fixed). We obtain that for large enough $x$ there are at least
\begin{equation*}
c_2(x/c_1)^{1/\delta} \ge x^A
\end{equation*}
variants for $T$, where the constant $A > 0$ is determined by $\delta$ and thus depends
only on the absolute type of $\G$.

It remains to recall that for each $T$ there exists a collection of parahoric subgroups
$P = (P_v)_{v\in V_f}$ such that $P_v = P_{0,v}$ for $v\in V_f\setminus T$,
$P_v \not\cong P_{0,v}$ for $v\in T$ and $P$ has $\cO$-maximal type. Each such collection
defines a maximal arithmetic subgroup of $H$ of covolume less than $x$ and
subgroups corresponding to different $T$'s are not conjugate. The number of maximal
arithmetic subgroups obtained this way is at least $x^A$ with $A>0$, a constant depending
only on the absolute type of $\G$. This proves the lower bound for
$m^u_H(x)$.
\subsection{} \begin{remark} Let us point out that
all the maximal arithmetic subgroups constructed in \ref{s51} are commensurable.
It is also possible to construct different commensurability classes which contain arithmetic
subgroups of covolumes less than~$x$. This may enlarge the constant $A$ in the asymptotic
inequality but, as it follows from the first part of the proof and the conjecture on the
number of isomorphism classes of fields with discriminant less than $x$, would hardly
change the type of the asymptotic.
\end{remark}
\subsection{Non-cocompact case.}\label{s52}
Let now $\Gamma_1$ be a non-uniform irreducible lattice in $H$ which exists by the assumption
of part~B of the theorem, and let $\G$ be a corresponding algebraic $k$-group. Arithmetic
subgroups of $H$ which are induced from $\G(k)$ will be all non-cocompact (they are all
commensurable with $\Gamma_1$). To provide a lower bound for $m^{nu}_H(x)$ it remains to
repeat the argument of~\ref{s51} for the group $\G$.

Note that in contrary to the cocompact case the existence of non-cocompact arithmetic subgroups
in $H$ is in general not guaranteed by the condition that $H$ is isotypic (a counterexample
was already mentioned in~\ref{s22}). Still the conditions under which such examples can be
constructed are rather exceptional and in most of the cases isotypic groups contain both cocompact
and non-cocompact arithmetic subgroups.

\bigskip

The theorem is now proven.

\section{Corollaries, conjectures, remarks}\setcounter{equation}{0}

\subsection{}\label{s61}
\begin{cor} There exists a constant $C_1$ which depends only
on the type of almost simple factors of $H$ such that if $\Lambda$ is a principal arithmetic
subgroup of $H$ and $\Gamma = N_H(\Lambda)$ has covolume less than $x$, then
$[\Gamma:\Lambda] \le x^{C_1}$.
\end{cor}

\begin{proof} By \cite{BP} (see Sect.~\ref{1:BP} for the notations and precise references):
\begin{equation}\label{e61}
[\Gamma : \Lambda] \le n^{\epsilon\# S}\cdot
2h_l^{\epsilon'} n^{\epsilon a(k) + \epsilon' a(l) +\epsilon \# T} (\D_l/\D_k^{[l:k]})^{\epsilon''}
\cdot \prod_{v\in V_f}\#\Xi_{\Theta_v}.
\end{equation}
Now, since $\mu(H/\Gamma)\le x$, the group $\Gamma$ has to satisfy the conditions
on the subgroups of covolume less than $x$ obtained in the proof of the upper
bound of Theorem~1 (in Theorem~1 only maximal arithmetic subgroups are considered
but the proof of the upper bound applies without a change to arbitrary principal
arithmetic subgroups and their normalizers thus providing a somewhat stronger result
to which we appeal here). We have:
\begin{gather*}
%\D_k \le (x/c_1)^{1/\delta_1},\ \
\D_l/\D_k^{[l:k]} = \D_{l/k} \le x^{c_1},\ \
a(k) \le c_2\log x {\rm\ \ (Sect.~\ref{s41})};\\
\# T \le \# T_1 \le c_{3}\log x {\rm\ (Sect.~\ref{s43})};\\
\{v\in V_f,\ \#\Xi_{\Theta_v}\neq 1\}\subset T_1\cup R\cup T_2 {\rm\ as\ for\ the\ rest\ of\ }
v,\ P_v {\rm \ is\ special,\ so}\\
\#\{v\in V_f,\ \#\Xi_{\Theta_v}\neq 1\} \le \#(T_1\cup R\cup T_2) \le c_4\log x
{\rm\ (Sect.~\ref{s44}).}
\end{gather*}
Also recall that $\#\Xi_{\Theta_v} \le r+1$, $r$ is the absolute rank of $\G$;
$h_l \le c^{[l:\Q]}\D_l\le x^{c_5}$ (see e.g. proof of Prop.~\ref{prop:BP}) and $a(l) \le
3a(k)$ (as $[l:k]\le 3$). Altogether these imply the proposition.
\end{proof}

\subsection{}\label{s62}
For some particular cases the bound in Corollary~\ref{s61} can be improved. Let
us assume that the degrees of the fields of definition of the arithmetic
subgroups are bounded:
\begin{equation}\label{e62}
[k:\Q] \le d,
\end{equation}
which is the case, for example, if we consider only non-uniform lattices in $H$.

Assumption (\ref{e62}) implies that the number $m$ of different prime
ideals $\cP_1, \ldots \cP_m$ of $\cO_k$ such that ${\rm Norm}
(\cP_1\ldots\cP_m) \le x$ is bounded by $c\log
x/\log\log x$, $c=c(d)$ (instead of the bound $\log x$ which we used for the
general case). Indeed, for $k = \Q$ it follows from the Prime Number Theorem,
and the case of arbitrary $k$ of bounded degree can be easily reduced to the
rational case.

Therefore assumption (\ref{e62}) implies that most of the terms in
(\ref{e61}) are \linebreak $\le c^{\log x/\log\log x}$ with $c=c(d)$. What
remains is $\D_{l/k} = \D_l/\D_k^{[l:k]}$ (for type $\Dn_r$, $r$ even) and
$h_l$ (or $h_k$ if $l=k$). The former, in fact, appears in the formula as an
upper bound for $2^{\# R}$ (see \cite[Sect.~5.5]{BP}), which again can be
improved to $c^{\log x/\log\log x}$ by the same argument. What remains is
the class number.

Going back to \cite[Sect.~5 and Prop.~0.12]{BP}, we see that what occurs in
the formula is not $h_l$ but the order of the group $C_n(l)$, which consists of
the elements of the class group $C(l)$ whose orders divide $n$ (as before
$n$ is a constant determined by the type of $H$). Instead of using the trivial
bound $\#C_n(l)\le h_l$, let us keep it as it is. We now come to the following
formula:
\begin{equation}\label{e63}
[\Gamma:\Lambda] \le c^{\log x/\log\log x} \#C_n(l),\ x\ge\mu(H/\Gamma).
\end{equation}
If $n = p$ is a prime, let $\rho_p(l)$ denote the $p$-rank of $C(l)$. Then,
clearly, $\#C_p(l) \le p^{\rho_p(l)}$, and in general for $n =
p_1^{\alpha_1}\ldots p_m^{\alpha_m}$, $\#C_n(l) \le p_1^{\alpha_1\rho_{p_1}(l)}
\ldots p_m^{\alpha_m\rho_{p_m}(l)}$. So we are interested in the upper bounds
for $p$-ranks of the class groups.

Apparently, even though this and related questions were much studied, there are
very few results beyond the celebrated Gauss' theorem which can be applied in our case.
We have:
\begin{itemize}
\item[$(i)$] if $[l:\Q] = 2$, then $\rho_2(l) \le t_l - 1$ (by Gauss);
\item[$(ii)$] if $[k:\Q] = 2$ and $[l:k] = 2$,  then $\rho_2(l) \le 2(t_l + t_k - 1)$
(by \cite[Theorem~2]{Cor});
\end{itemize}
where $t_k$ (resp. $t_l$) denotes the number of primes ramified in $k/\Q$ (resp. $l/k$).

From this we obtain
\begin{equation}\label{e64}
\#C_n(l) \le n^{c \log x/\log\log x},
\end{equation}
if $n$ is a power of $2$ and $l$ is as in $(i)$ or $(ii)$.

Similar results for other $n$ and other fields can be only conjectured: Even
$\rho_3(k)$ for quadratic fields $k$ seems to be out of reach with the currently
available methods. Still estimates (\ref{e63}) and (\ref{e64}) imply the
following:

\subsection{} \label{s63}
\begin{cor} Let $H$ be a simple Lie group of type $\An_{2^\alpha-1}$ $(\alpha > 1)$,
$\Bn_r$, $\Cn_r$, $\Dn_r$ $( r\neq 4 )$, $\En_7$, $\En_8$, $\Fn$ or $\Gn$.
There exists a constant $C_2$ which depends only on the type of $H$ such that
if $\Lambda$ is a non-cocompact principal arithmetic subgroup of $H$ and
$\Gamma = N_H(\Lambda)$ has covolume less than $x$, then
$[\Gamma:\Lambda] \le C_2^{\log x/\log\log x}$.
\end{cor}

\begin{proof} Indeed, the assumption that $H$ has one of the given types
implies that $n$ is a power of $2$ (see the definition of $n$ in Sect.~\ref{1:BP}).
Since $H$ is simple and the arithmetic subgroup is non-compact, its field of
definition $k$ is either $\Q$ or an imaginary quadratic extension of $\Q$,
depending on $H$ being real or complex Lie group (see also Sect.~\ref{s42}). Finally,
the type of $H$ is not $\Dn_4$ implies $[l:k]\le 2$. The corollary now follows
from the discussion in~\ref{s62}.
\end{proof}

We expect similar estimates to be valid for the non-uniform lattices in other
groups but we do not know how to prove it.

\subsection{} \begin{remark}
Concerning the general case, let us point out that if the degrees
of the fields are a priori not bounded then we can not expect a $\log
x/\log\log x$--bound for the $p$-rank of the class group. An example of a
sequence of fields $k_i$ for which $\rho_2(k_i)$ grows as $\log \D_{k_i}$ was
constructed by F.~Hajir (On the growth of $p$-class groups in $p$-class field
towers, {\em J.~Algebra} {\bf 188} (1997), 256--271), the fields $k_i$ in
Hajir's example form an infinite class field tower. This remark together with
the previous estimates motivates the following question:

{\it Is the estimate in Corollary~\ref{s61} sharp, i.e. given a group $H$ is
there a constant $C_0 = C_0(H) > 0$ such that there exists an infinite sequence
of pairwise non-conjugate principal arithmetic subgroups $\Lambda_i$ in $H$ for
which $[\Gamma_i:\Lambda_i] \ge \mu(H/\Gamma_i)^{C_0}$, where $\Gamma_i =
N_H(\Lambda_i)$ and $\mu$ is a Haar measure on $H$?}
\end{remark}

Corollaries~\ref{s61} and \ref{s63} are important for~\cite{BL} where we
study the growth rate of the number of irreducible lattices in semi-simple Lie
groups.

\subsection{} \label{s64} \begin{remark}
All along the line the groups of type $\An_1$ were excluded. The reason for this is that we
can hardly hope to use the formula from Sect.~\ref{2:BP} combined with an analogue of
Proposition~\ref{prop:BP} for this case even to prove a finiteness result. Still one can follow another
method, also due to Borel, and use geometric bounds for the
index of a principal arithmetic subgroup in a maximal arithmetic. This indeed allows to establish the
finiteness of the number of arithmetic subgroups of bounded covolume in $\SL(2,\R)^a\times\SL(2,\C)^b$
\cite{B}. Now the problem is that the quantitative bounds which can be obtained this way are only
exponential. We suppose that the true bounds should be similar to the general case (and conjecturally
polynomial), although we do not know how to prove this conjecture and leave it as an open problem:

{\it Find the growth rate of the number of maximal arithmetic subgroups for the semi-simple Lie groups
whose almost simple factors have type $\An_1$ or obtain a better than exponential upper bound for the
growth.}
\end{remark}

\subsection{}\label{s65}
Two conjectures were mentioned in the introduction:
\begin{conjecture}
There exists an absolute constant $B$ such that for large enough $x$ the number of isomorphism classes
of number fields with discriminants less than $x$ is at most $x^B$.
\end{conjecture}
\begin{conjecture}
Given a connected semi-simple Lie group $H$ without almost simple factors of type $\An_1$ and compact
factors, there exists a constant $B_H > 0$ which depends only on the type of almost simple factors of
$H$ such that for large enough $x$ the number of conjugacy classes of maximal irreducible arithmetic
subgroups of $H$ of covolumes less than $x$ is at most $x^{B_H}$.
\end{conjecture}
We now prove:

\begin{prop}
Conjectures $1$ and $2$ are equivalent.
\end{prop}
\begin{proof} $'1\to 2'$ follows directly from the proof of the upper bounds in Theorem~1.

\noindent
$'2\to 1'.$ Assume that Conjecture~2 is true but Conjecture~1 is false, i.e. $m^u_H(x) + m^{nu}_H(x)
\le x^{B_H}$ for every $x > x_0$, and for an arbitrary $C$ there exists $x>x_0$ such that $N(x) > x^C$.
We shall need some additional assumption on $x_0$ which will become clear later but could be imposed
from the beginning. So let us fix $C>1$ and let $x > x_0$ be such that $N(x) > x^C$.

Let $N_{i,j}(x)$ denote the number of extensions of $\Q$ of discriminant less than $x$ which have
precisely $i$ real and $j$ complex places. We have
\begin{equation*}
N(x) = \sum_{\substack{i=1,\ldots,n\\j=1,\ldots,m}}N_{i,j}(x).
\end{equation*}
The condition that the discriminants of the fields are less than $x$ implies by
Minkow\-ski's theorem that the degrees of the extensions are bounded by $c\log
x$ for an absolute constant $c$, and so the number of summands is less than
$(c\log x)^2$. By Dirichlet's box principle there exists a pair $(i,j)$ such
that $N_{i,j}(x) > x^C/(c\log x)^2 \ge x^{C-1}$ (this inequality requires
$(c\log x)^2 \le x$ which is true for large enough x and gives a first
condition on $x_0$). Let $\mathcal K$ be the set of such number fields,
$\#\mathcal K = N_{i,j}(x) > x^{C-1}$. Consider a simply connected semi-simple
Lie group $H$ which has $i$ split real simple factors and $j$ complex simple
factors all of the same type. For each $k\in\mathcal K$, let $\G/k$ be a simply
connected, absolutely simple split group of the same absolute type as the
simple factors of $H$ and defined over $k$. Let $P = (P_v)_{v\in V_f}$ be a
coherent collection of parahoric subgroups of $\G$ all of which are
hyperspecial (such a collection exists since $\G$ splits over $k$), and let
$\Lambda$ be the principal arithmetic subgroup of $H$ defined by $P$. We have:
\begin{itemize}
\item[a)] for $S=V_\infty(k)$, $\G_S\cong H$;
\item[b)] $\mu(H/\Lambda) = \D_k^{\d}\left(\prod_{i=1}^{r}\frac{m_i!}{(2\pi)^{m_i+1}}\right)^{[k:\Q]}\E(P)$
\ by Prasad's formula.
\end{itemize}

Using the orders of finite groups of Lie type the Euler product $\E(P)$ can be
expressed as a product of the Dedekind zeta function of $k$ and certain Dirichlet $\rm L$-functions
at the integers $m_i+1$, $m_i$ are the Lie exponents of $\G$ (see \cite[Rem.~3.11]{P}).
Obvious inequalities $L(s,\chi)\le\zeta_k(s)$ and $\zeta_k(s)\le \zeta(s)^{[k:\Q]}$, for $s \ge 2$, imply
that there exists a constant $c_1$ which depends only on the type of simple factors of $H$ and such that
each zeta or $\rm L$-function in the product is bounded from above by $c_1^{[k:\Q]}$.
Since $[k:\Q]\le c\log x$, we have $\E(P) \le (c_1^{c\log x})^r$. By definition of the
set $\mathcal K$, $k\in \mathcal K$ implies $\D_k\le x$. Therefore we obtain
\begin{equation*}
\mu(H/\Lambda) \le x^{\d} c_2^{c\log x} (c_1^{c\log x})^r \le x^\delta,
\end{equation*}
where $\delta$ is $>1$ and depends only on the type of simple factors of $H$. The latter inequality may
require that x is larger than certain value which depends on the type of simple factors of $H$ and gives us
the second condition on $x_0$. Clearly, both conditions do not depend on $C$ and so could be imposed from
the beginning.

For each $k\in\mathcal K$ we have at least one maximal arithmetic subgroup of $H$ of covolume less then
$x^\delta$. Now if we take $C = \delta B_H + 1$ we arrive to a contradiction with Conjecture~2 for $H$
and $x^\delta > x_0$.
\end{proof}
Let us remark that in the proof Conjecture~2 is used only for non-cocompact arithmetic subgroups
of semi-simple groups $H$ which have simple factors of a fixed type, it then implies Conjecture~1 which
in turn implies Conjecture~2 in the whole generality. It is possible to specify further the relation
between two conjectures but we shall not go into details. What we would like to emphasize is that our
result provides a new geometric interpretation for a classical number theoretic problem. An optimistic
expectation would be that study of the distributions of lattices in semi-simple Lie groups
can give a new insight on the number fields and their discriminants.

\appendix
\section*{Appendix}
\begin{center} \it by J. Ellenberg and A. Venkatesh \end{center}
\medskip
\stepcounter{section}

\subsection{}\label{thmone}
Let $N(X)$ denote the number of isomorphism
classes of number fields with discriminant less than $X$.
%\medskip

\begin{thm}
For every $\epsilon > 0$ there
is a constant $C(\epsilon)$ such that
$\log N(X) \leq C(\epsilon) (\log X)^{1+\epsilon}$,
for every $X \geq 2$.
\end{thm}

In fact we prove the more precise upper bound that
$$%\beq
\log N(X) \leq C_6 \log X \exp(C_7 \sqrt{\log \log X})
$$%\eeq
for absolute constants $C_6,C_7$.

This theorem (almost) follows from \cite[Theorem 1.1]{EV},
the only point being to control the dependence
of implicit constants on the degree of the number field.

We refer to \cite{EV} for further information and for some motivational
comments about the method.  In the proof $C_1, C_2, \dots$ will denote certain
{\em absolute} constants.

\subsection{}
Let $K$ be an extension of $\Q$ of degree $d \geq 200$.
% Let $\order_K \subset K$ be the ring of integers.
Denote by $\Sigma(K)$ the set of embeddings of $K$ into  $\C$ ($\#\Sigma(K) = d$),
and by $\overline{\Sigma}(K)$ a set of representatives for $\Sigma(K)$ modulo complex
conjugation (in the notations of the paper $\overline{\Sigma}(K) = V_\infty (K)$).
% Thus $\#\Sigma(K) = d$.
% For $\sigma \in \overline{\Sigma}(K)$ we put $K_{\sigma} = \R$ or $\C$ according to whether
% $\sigma$ is real or complex.
We regard the ring of integers $\order_K$ as a lattice
in $K \otimes_{\Q} \R = \prod_{\sigma \in \Sigma(K)} K_{\sigma}$. We endow the real vector
space $K \otimes_{\Q} \R$
with the supremum norm, i.e. $\|(x_\sigma)\| = \sup_{\sigma}
|x_{\sigma}|$.
Here $|\cdot|$ denotes the standard absolute value on $\C$.
In particular, we obtain a ``norm''
on $\order_K$ by restriction. Explicitly, for $z \in \order_K$, we have
$\|z\| = \sup_{\sigma \in \Sigma(K)} |\sigma(z)|$.

We denote by $\M_d(\Z)$ (resp. $\M_d(\Q)$) the algebra
of $d$ by $d$ matrices over $\Z$ (resp. $\Q$).

By the ``trace form'' we mean the pairing $(x,y) \mapsto \Tr(xy)$.
It is a symmetric nondegenerate $\Q$-bilinear pairing on $K^2$.

Let $s$ be a positive integer which we shall specify later.
We denote by \linebreak
$\y = (y_1, y_2,\ldots,y_s)$ an ordered $s$-tuple of elements of $\order_K$
and write $\|\y\| :=$ \linebreak
$\max(\|y_1\|,\ldots, \|y_s\|)$. For $\y = (y_1, \ldots, y_s) \in \order_K^s$
and $l \geq 1$, we shall set
\begin{gather*} %equation} \begin{aligned}
S(l) = \{(k_1,\ldots,k_s) \in \Z^s:k_1 + \ldots+ k_s \leq l,\ k_1, \ldots, k_s \geq 0\},  \\
S(\y, l) = \{y_1^{k_1} y_2^{k_2} \ldots y_s^{k_s} : (k_1, \ldots, k_s) \in S(l)\} \subset \order_K.
\end{gather*} %aligned} \end{equation}
If $S$ is a subset of $S(l)$ we denote by $S(\y)$ the set $\{y_1^{k_1}
y_2^{k_2} \ldots y_s^{k_s} : (k_1, \ldots, k_s) \in S$.\}

\subsection{} \begin{lemma} \label{lem:span}
Let $S$ be a subset of $S(l)$ such that $S(\y)$ spans a $\Q$-linear
subspace of $K$ with dimension strictly greater than $d/2$.  Let $S+S$
be the set of sums of two elements of $S$.  Then $(S+S)(\y)$ spans $K$
over $\Q$.
\end{lemma}
\proof (cf. \cite[Lemma 2.1]{EV}.)
Suppose that there existed $z \in K$ which was perpendicular, w.r.t.
the trace form, to the $\Q$-span of $(S+S)(\y).$
Since $(S+S)(\y)$ consists precisely of all products
$\alpha \beta$, with $\alpha, \beta \in S(\y)$,
it follows that
\begin{equation} \label{eqn:bart}\Tr(z \alpha \beta) = 0, \ \ (\alpha, \beta
\in S(\y)).\end{equation}
Call $W \subset K$ the $\Q$-linear span of $S(\y)$.
Then (\ref{eqn:bart}) implies that $z W$ is perpendicular to $W$
w.r.t. the trace form, contradicting $\dim(W) > d/2$. \qed

\subsection{} \begin{lemma} \label{lem:traces}
Let $\mathcal{C} \subset \order_K$
be a finite subset containing $1$ and generating $K$ as
a field over $\Q$. Let
$z_1, z_2, \dots, z_d$ be a $\Q$-linear basis for $K$.
For each $u \in \mathcal{C}$,
let $ M(u) =
\left(\mathrm{Tr}_{K/\Q}(uz_i z_j)\right)_{1 \leq i,j\leq d}
\in \M_d(\mathbb{Q}).$
Then the $\Q$-subalgebra of $\M_d(\Q)$
generated by $M(u) M(1)^{-1}$,
as $u$ ranges over $\mathcal{C}$, is isomorphic to $K$.
\end{lemma}
\proof
(cf. \cite[Lemma 2.2]{EV}.) In fact, $M(u) M(1)^{-1}$ gives the matrix
of ``multiplication by $u$,'' in the basis $\{z_i\}$. \qed

\subsection{}
\begin{lemma} \label{lem:redtheory}
There is an absolute constant $C_1 \in \mathbb{R}$
such that,
for any $K$ as above,
there exists a basis $\gamma_1, \gamma_2, \dots, \gamma_d$
for $\order_K$ over $\Z$ such that
\begin{equation} \label{eqn:red} \|\gamma_j\| \leq
\|\gamma_{j+1}\|, \ \
\prod_{i=1}^{d}
\|\gamma_d\| \leq \disc^{1/2}  C_1^{d}, \ \
\|\gamma_i\| \leq (C_1^d \disc^{1/2})^{\frac{1}{d-i}}  \ (i<d).
\end{equation}
\end{lemma}
($\disc$ denotes the absolute value of the discriminant of $K$.)
\proof This is reduction theory (cf. \cite[Prop.~2.5]{EV}). The final
statement of (\ref{eqn:red}) follows from the preceding
statements, in view
of the fact that $\|\gamma_j\| \geq 1$ for each $j$.
\qed

\subsection{}
Let $r,l$ be integers such that $d/2 < r \leq |S(l)| = \binom{l+s}{s}$.

\begin{lemma} \label{lem:dimcount}
Suppose $W \subset K$ is a $\Q$-linear subspace
of dimension $r$, and let $S \subset S(l)$ be a subset of size
$r$.  Then there exists $\y= (y_1, y_2, \ldots y_s) \in W^s$
such that the elements of $S(\y)$ are $\Q$-linearly independent.
\end{lemma}

\proof
This is precisely \cite[Lemma 2.3]{EV}.
\qed

\subsection{} \begin{lemma} \label{lem:lambda}
Let $\Lambda = \Z \gamma_1 + \Z \gamma_2 + \dots + \Z \gamma_{r}$ and
let $S \subset S(l)$ be a subset of size $r$.
Then there is $\y = (y_1, y_2, \ldots, y_s) \in \Lambda^s$
such that the elements of $S(\y)$ are linearly independent over $\Q$,
and $\|\y\| \leq r^2 l (C_1^d \disc^{1/2})^{\frac{1}{d-r}}$.
\end{lemma}
\proof
Considering $\Lambda^s$ as a $\Z$-module of rank $rs$, the proof of
\cite[Lemma 2.3]{EV} shows that there is a polynomial $F$ of
degree at most $rl$ in the $rs$ variables so that the elements of
$S(\y)$ are linearly independent over $\Q$ whenever $F(\y) \neq 0$.
Lemma 2.4 of \cite{EV} then shows that we can choose such a $\y$ whose
coefficients are at most $(1/2)(rl+1) \leq rl$.  It follows that
\beq
\|y_i\| \leq r^2 l (C_1^d \disc^{1/2})^{\frac{1}{d-r}}
\eeq
for $i = 1,2,\ldots s$.
\qed

\subsection{} \begin{lemma} \label{lem:final}
The number of number fields with degree $d \geq 200$ and discriminant
of absolute value at most $X$ is at most
\beq
(C_3 d)^{d \exp(C_4 \sqrt{\log d})} X^{\exp(C_5 \sqrt{\log d})}.
\eeq
\end{lemma}
\proof
Fix once and for all a total ordering of $S(2l)$. We denote
the order relation as $(k_1,\ldots,k_s) \prec (k_1', \ldots, k_s')$.
Choose $S \subset S(l)$ of cardinality $r$ as above.

Let $K$ have degree $d$ over $\Q$
and satisfy $\disc < X$. Chose $\y$ as in Lemma~\ref{lem:lambda}.
By Lemma~\ref{lem:span}, $S(2l)(\y)$
spans $K$ over $\Q$. It follows that there exists a subset $\Pi
\subset S(2l)$ of size $d$
such that $\{z_1, \dots, z_d\} := \{y_1^{k_1} y_2^{k_2} \ldots y_s^{k_s}: (k_1, k_2, \ldots,
k_s) \in \Pi\}$ forms a $\Q$-basis for $K$, and such that the ordering
$z_1,\ldots, z_d$ conforms with the specified ordering on $\Pi \subset
S(2l)$.

We apply Lemma~\ref{lem:traces} to
$\{z_1, \dots, z_d\}$ and $\mathcal{C} = (1, y_1, y_2, \ldots, y_s)$.
Then each product $u z_i z_j \, (u \in \mathcal{C}, 1 \leq i,j \leq d)$
is contained in $S(4l+1)$.

Put $\bfA = (\Tr(y_1^{k_1} y_2^{k_2} \ldots y_s^{k_s}))_{(k_1, k_2,
  \ldots k_s) \in S(4l+1)}.$  For each $K$, the collection of matrices $M(u)$ is determined by
$\bfA$ and $\Pi$.
Since $|\Tr(z)| \leq d \|\y\|^{4l+1}$ for any
$z \in S(\y,4l+1)$,
the number of possibilities for $\bfA$ is
at most $(d \|\y\|^{4l+1})^{|S(4l+1)|}$;
since $\Pi$ is a subset of $|S(2l)|$,
the number of possibilities for $\Pi$ is at most
$2^{|S(2l)|}$.

Lemma~\ref{lem:traces} now yields that the number of possibilities for
the isomorphism class of $K$ is at most $2^{|S(2l)|} (d
\|\y\|^{4l+1})^{|S(4l+1)|}$.  By our bound on $\|\y\|$ we now have that the
number of possibilities for $K$ is at most
\begin{equation}
2^{|S(2l)|} (d (r^2 l (C_1^d \disc^{1/2})^{\frac{1}{d-r}})^{4l + 1})^{|S(4l+1)|}.
\label{eq:bigineq}
\end{equation}
Note that $|S(4l+1)| = \binom{s + 4l+1}{s}$.

Now, just as in the paragraph following (2.6) of \cite{EV}, we choose
$s$ to be the greatest integer less than $\sqrt{\log{d}}$ and $l$ to
be the least integer greater than $(ds!)^{1/s}$.  Note that $l <
\exp(C_2 \sqrt{\log d})$.  Now
$|S(l)| = \binom{s+l}{s}$ is at least $d$, so we may choose $r$
between $d/2$ and $3d/4$.  In particular, $r^2 l < d^3.$  Also, $\binom{s+4l+1}{s}$ is at most $10^s d$
and $|S(2l)| = \binom{s+2l}{s} \leq 6^s d$.  Finally, $s < 2
\sqrt{\log d}$.

Substituting these values into \eqref{eq:bigineq} we get that the
number of possible $K$ is at most
\beq
2^{6^s d} (d (d^3 (C_1^d X^{1/2})^{4/d})^{5 \exp(C_2 \sqrt{\log d}})^{10^s d}
\eeq
which is in turn at most
\beq
(C_3 d)^{d \exp(C_4 \sqrt{\log d})} X^{\exp(C_5 \sqrt{\log d})}.
\eeq
\qed

\subsection{} \begin{prop}
There are absolute constants $C_6,C_7$ with
\beq
\log N(X) \leq C_6 \log X \exp(C_7 \sqrt{\log \log X}).
\eeq
\end{prop}
\proof
By Minkowski's discriminant bound, there is an absolute constant $C_6 > 1$
such that $\disc > C_6^{[K:\Q]}$ for any extension $K/\Q$; so we may
take $d$ to be bounded by a constant multiple of $\log X$.
From Lemma~\ref{lem:final}
it now follows that the logarithm of the  number of extensions $K/\Q$
with $\disc < X$ and $[K:\Q] \geq 200$
is bounded by $C_6 \log X \exp(C_7 \sqrt{\log \log X})$.
Trivial bounds suffice to show that
the number of $K$ with $\disc < X$ and $[K:\Q] < 200$
is $\leq C_8 X^{200}$.
\qed

\end{document}